\documentclass[12pt,reqno]{amsart}

\usepackage[a4paper,margin=1in]{geometry}
\usepackage{amsfonts,amssymb,amsmath,amsthm}
\usepackage{graphicx}
\usepackage{tikz}
\usepackage{psfrag}
\usepackage{calc}
\usepackage{epic}
\usepackage{hyperref}
\usepackage{float}
\usepackage{latexsym}
\usepackage{color}
\usepackage{subfigure,dsfont}
\usepackage{todonotes}

\newtheorem{theorem}{Theorem}
\newtheorem{lemma}[theorem]{Lemma}
\newtheorem{proposition}[theorem]{Proposition}
\newtheorem{corollary}[theorem]{Corollary}
\newtheorem{definition}[theorem]{Definition}

\theoremstyle{definition}
\newtheorem{example}{Example}
\newtheorem{remark}{Remark}

\newcommand{\be}{ \begin{equation}}
\newcommand{\ee}{\end{equation}}
\newcommand{\ben}{ \begin{equation*}}
\newcommand{\een}{\end{equation*}}

\newcommand{\ba}{ \begin{align}}
\newcommand{\ea}{\end{align}}
\newcommand{\ban}{ \begin{align*}}
\newcommand{\ean}{\end{align*}}

\newcommand{\Prob}{\mathbb P}

\newcommand{\E}{\mathbb E}

\newcommand{\R}{\mathbb{R}}

\newcommand{\I}[1]{{\bf 1}{\left\{#1\right\}}} 

\def\R{\mathbb{R}}

\def\L{\mathcal{L}} 

\let\oldsum\sum
\renewcommand{\sum}{\displaystyle\oldsum} 

\title{A Distributional Equation Approach to Front Propagation}

\author{Matthieu Jonckheere}\address{LAAS-CNRS, Université de Toulouse, CNRS, Toulouse, France}\email{matthieu.jonckheere@laas.fr}
\author{Seva Shneer}\address{School of MACS, Heriot-Watt University, Edinburgh, EH14 4AS, United Kingdom}\email{V.Shneer@hw.ac.uk}

\begin{document}

\begin{abstract}
Travelling waves of mean-field particle systems are usually studied through
the partial differential equation satisfied by the mean-field limit. We
develop a different and entirely probabilistic set of tools. The wave profile
is characterised as the solution of a distributional fixed-point equation for
a tagged particle seen in the frame of the wave, and that equation is solved by
relating it to a linear one --- a smoothing transform --- whose solutions are
the martingale limits of a branching random walk. 
We apply this approach to a range of models, recover some known results and obtain a number of new ones. We first study models where particles evolve according to independent L\'evy processes and, at Poisson epochs synchronise randomly selected pairs of particles (the lower of the two moves instantaneously to the position of the higher of the two). 
We characterise the speeds and travelling waves for all such models. In particular, this covers the case where the L\'evy process is a continuous-time random walk with an arbitrary jump distribution; previously, only exponentially distributed jumps had been treated. We also consider a variant in which, instead of synchronising, the lower of the two particles jumps up by a random amount, drawn independently at each jump.
This model is motivated by the so-called power-of-2 systems from queueing theory. Finally, we consider a new model where updates are given by the mixture of synchronisation and power-of-2: each update is a random value but the resulting position of the lower particle is capped by that of the higher one. 
\\
In all cases, we characterise the speed, the travelling waves and their tail behaviour; the probabilistic nature of our approach makes the tail analysis particularly tractable.

\end{abstract}

\maketitle

\section{Introduction}

Reaction-diffusion systems with mean-field interactions arise in many areas. This work is motivated by their appearance in statistical
physics, where local interactions between particles propagate into macroscopic
patterns \cite{BD2,Liggett}; in population dynamics, where they describe the
spatial spread of populations subject to reproduction, competition and
migration \cite{Balazs2014,Greenberg1995}; and in queueing theory and
communication networks, where interacting servers affect each other's
workloads through load-balancing mechanisms \cite{Stolyar2023TW}. In the limit
of infinitely many particles the microscopic interactions collapse into a
deterministic macroscopic description, whose most striking emergent feature is
the formation of travelling waves: self-sustaining, shape-preserving fronts
moving at constant velocity. Travelling waves provide the link between the
microscopic stochastic dynamics and the macroscopic equation
\cite{BD,BD2,GJ2018}, and raise the classical questions of speed selection,
stability of the profile, and existence of such coherent structures.

More precisely, the models we consider are of the following type. Particles move according to a L\'evy
process and interact through one of three mechanisms: \emph{synchronisation},
where the lower of two interacting particles jumps to the position of the
higher, which in the diffusive case produces the F-KPP equation; the
\emph{power-of-2} growth mechanism, where a particle which finds itself below
a uniformly sampled competitor is pushed up by an independent random amount; and finally a \emph{mixed} mechanism where the power-of-2 update is capped at the level of the competitor. We will provide exact definitions of all the models in what follows.

The usual route to the establishment of the existence and uniqueness of travelling-wave solutions and to the understanding of their properties is the partial differential equation
satisfied by the mean-field limit. We propose a different one. Write $\chi$ for the
position of a \emph{tagged particle} of the mean-field (McKean--Vlasov) dynamics,
seen in the frame moving with the wave, and follow it from one of its interactions to the
next: what results is a closed distributional equation for the law of $\chi$,
and that law is precisely the travelling-wave profile. The equation is an instance of a random fixed-point one, and we discuss strategies for constructing its solutions and studying their properties.

It is not our goal to demonstrate that our approach is a substitute for the PDE one but rather to propose alternative, probabilistic way of constructing travelling-wave solutions. In addition, our approach leads to new findings: the travelling waves can be constructed using martingale limits of associated branching processes and their tails can be characterised accordingly.

\subsection*{The two distributional equations}

Two families of random fixed-point equations are crucial for our analysis. The first,
which the tagged particle produces directly, is
\begin{equation} \label{eq:general_type_max}
\chi \overset{d}{=} \max\{\chi_1,\chi_2\} - B,
\end{equation}
where $\chi_1,\chi_2$ are independent copies of $\chi$ and $B$ is independent of
them. The second is
\begin{equation} \label{eq:general_type_sum}
\chi \overset{d}{=} (\chi_1+\chi_2)A,
\end{equation}
again with $\chi_1,\chi_2$ independent copies of $\chi$ and all the variables on
the right-hand side independent, and now non-negative. We refer to
\eqref{eq:general_type_max} and \eqref{eq:general_type_sum} as the non-linear
and the linear equation respectively.

{It is easy to see how the approach identified above leads to an equation of type \eqref{eq:general_type_max} in the case of synchronisation. In this model
$B=cT-\Gamma(T)$, with $c$ denoting the speed of the wave, $T$ denoting the Exponential time between the two interactions involving our tagged particle and $\Gamma(T)$ denoting the displacement of the particle by a L\'evy process in time $T$. It turns out that the power-of-2 mechanism gives, after the change of variable
$\eta=-\chi$ and some further manipulations, an equation of the same type, with $B$ built from the
equilibrium jump size and the driver's all-time supremum in the wave frame. The model of Section \ref{sec:mixed} produces an equation which is a slight modification of \eqref{eq:general_type_max}.
Different mechanisms thus lead to distributional equations of the same type, and our approach allows therefore to treat these models in a similar way.}

Equations of type \eqref{eq:general_type_max} have received limited attention
in the literature. Equations of type \eqref{eq:general_type_sum} --- the
smoothing transform --- have on the contrary been studied extensively
\cite{Liu2001,liu1999asymptotic,Bigginsetal91}, and their solutions are known
to be martingale limits of branching random walks. The bridge between the two
is Lemma~\ref{lemma:deq_main}: a random variable $\chi$ with $\Prob(\chi\le x)=\E\bigl[e^{-Ve^{-\gamma x}}\bigr]$ solves
\eqref{eq:general_type_max}
if and only if $V$ solves \eqref{eq:general_type_sum} with $A=e^{-\gamma B}$.

\subsection*{Main results}

Our starting
point, Proposition~\ref{prop:TW-FP}, applies to any interaction which displaces the
particle initiating it and is translation covariant. It states that a
probability measure $h$ is a travelling wave of the mean-field equation with
speed $c$ if and only if
\[
\chi\;\overset{d}{=}\;\chi_1+\bar\theta(\chi_1,\chi_2)+\Gamma(T)-cT ,
\]
where $\chi,\chi_1,\chi_2$ are independent with law $h$ ($\Prob(\chi\le x)=h(x)$), $T\sim\mathrm{Exp}(1)$ is the
time to the tagged particle's next interaction, $\Gamma(T)$ is the driver's
displacement over that time, and $\bar\theta(x,y)$ is the displacement of a
particle at $x$ interacting with a partner at $y$. The equation expresses the
stationarity of a tagged particle in the frame of the wave, and its proof
requires no assumption on the underlying particle system. The three mechanisms of this
paper are three choices of $\bar\theta$, and everything that follows is the
analysis of this one equation.

The connection with the particle systems is made by
Theorem~\ref{prop:tw}. For finitely many particles the system has no fixed
origin, so one works with the configuration seen from its barycentre, which
has a random evolution around a linear drift. Under stationarity of this
recentred system, propagation of chaos in the stationary regime and a law of
large numbers for the barycentre, the stationary profiles converge to a
solution of the equation above, hence to a travelling wave: the one selected
by the finite systems.


Lemma~\ref{lem:bary-gen} shows that the barycentre law of large numbers follows
from three ingredients: stationarity, two-particle propagation of chaos and second-moment bounds on the profile
and on the interaction kernel, for an arbitrary mechanism and an arbitrary
L\'evy driver. 

The speed can then be identified intrinsically as
$c=\gamma_1+\bar\beta$ with $\gamma_1 = \E(\Gamma(1))$ and $\bar\beta=\E_{h^{\otimes2}}\bar\theta(\chi_1,\chi_2)$
the mean interaction displacement. 

\

We then demonstrate that the approach proposed here allows one to derive the well-known shape and speed of the F-KPP equation which one obtains in synchronisation model where the driver is a driftless Brownian motion.
The method also allowed us to obtain new results for other models with synchronisation. For instance, for models combining independent jumps at rate
$\lambda$ (in other words, independent continuous-time random walks) with synchronisation, we characterise the travelling wave for every
jump law possessing an exponential moment, with minimal speed
$c^*=\min_{\gamma}\gamma^{-1}\bigl(1+\lambda(\E[e^{\gamma X}]-1)\bigr)$ and an
explicit branching representation of the profile at every admissible speed
(Proposition~\ref{thm:BS_TW}). Previously the wave was known only for exponentially
distributed jumps \cite{baryshnikov2025large}.
Our method is very versatile and we note that it can easily be generalised to handle synchronisation of an arbitrary number of particles, or even a random number of particles picked independently from the same distribution at every interaction instance.

\

We then apply the same ideas to the power-of-2 dynamics. First, Proposition~\ref{prop:equiv_power2} states that, for an arbitrary driving
L\'evy process $\Gamma$ with $\E|\Gamma(1)|<\infty$, the power-of-2 profile
satisfies
\begin{equation}\label{eq:intro-decomp}
\chi-\overline H\;\overset{d}{=}\;\min\{\chi_1,\chi_2\}+\overline X
+\sup_{t\ge0}\bigl(\Gamma(t)-ct\bigr),
\end{equation}
all four variables independent, where $\overline X$ is the equilibrium law of
the jump size and $\overline H$ that of the descending ladder height of
$\Gamma(t)-ct$. 
When $\Gamma$ has no negative jumps {(which is a natural assumption in the queueing context)},
$\overline H\equiv0$ and the supremum has the Pollaczek--Khinchine law (for
driftless Brownian motion, for instance, it is exponential). {In this case \eqref{eq:intro-decomp} is also of the form
\eqref{eq:general_type_max} after the change of variable $\eta=-\chi$, which is
the structural connection between the power-of-2 and synchronisation models.}


{One can identify the speed $c$ to be equal to $\E\,\Gamma(1)+\E[X]/2$, which is intuitively clear. What our method allows to find however is the characterisation of the law of $\chi$, which is the travelling-wave profile. This has not been known before.}

{The probabilistic approach proposed in this paper allows to shed light on the tails of the travelling-wave profile}. The profile
is a Gumbel mixture $1-\E[e^{-We^{-\gamma x}}]$, and one tail decays at
rate $\gamma$. The other {one however} does not. It is governed by $\mu^*= \sup\{\mu : \E[e^{\mu A}]<\infty \}$ {with $A$ defined in \eqref{eq:general_type_sum}}.

For power-of-2 the inequality reverses: in the
explicit case of exponential jumps and unit Brownian variance the profile
decays like $x\,e^{-x}$ (Example~\ref{ex:pow2-tail}); the factor $x$
comes from the law of $A$, not from criticality, since the power-of-2 wave is
always supercritical.
{We also note that the method does have limitations.} When the driver has negative jumps, \eqref{eq:intro-decomp} still holds but the
extra factor $\overline H$ carries the fixed point outside the class of
smoothing transforms, and forces a strictly heavier tail than the one-sided
formula predicts (Remark~\ref{rem:levy_case}). {The exact shape of the tail in that case is stated as an open
problem.}

Section~\ref{sec:mixed} treats a model in which a single update mixes the two
mechanisms: a particle jumps, but not past a uniformly sampled competitor.
{This is a new model. Rather surprisingly, our approach allows to find a closed-form expression for the speed
\[
c=\frac{\lambda_2}{\log\bigl(1+\lambda_2/\lambda_1\bigr)} ,
\]
in terms of the two update rates $\lambda_1$ and $\lambda_2$ alone and (even more surprisingly) independently of the jump law. We also characterise the travelling-wave solution as a max-type equation with a
\emph{random} number of branches (i.e. satisfying a minor generalisation of equation \eqref{eq:general_type_max}),  that number being Poisson conditioned to be
positive (Proposition~\ref{prop:mixed})}.

\subsection*{Organisation}
Section~\ref{sec:mf} sets up the mean-field framework and the examples of interaction,
characterises travelling waves through the tagged-particle equation, and
connects this equation to the stationary regime of the particle systems. {We also discuss in detail the relevant literature and previously obtained results}. 
Section~\ref{sec:distr_equations} introduces the
two distributional equations and the correspondence between them,
Section~\ref{sec:prob-FKPP} treats the synchronisation mechanism,
Section~\ref{sec:power-of-2} the power-of-2 mechanism and
Section~\ref{sec:mixed} a model combining the two. Finally
Appendix~\ref{app:ladder} collects the input we need from fluctuation theory.

\section{Mean-field limits of reaction-diffusion systems and travelling waves}\label{sec:mf}

\subsection{Mean-field limits}

We consider a generic reaction-diffusion particle system where the diffusion dynamics are driven by a Lévy process $(\Gamma_t)_{t \ge 0}$ and the
reaction mechanism takes place at exponential times. More precisely, given $\L_d$ the infinitesimal generator of a Lévy process, let us define $(\xi^N(t))_{t\ge0}$ a Markov process in the state space $\R^N$ with generator given by
\begin{equation}
\label{eq:generator}
\L f(\xi)=\L_d f(\xi) + \frac{1}{ (N-1)} \sum_{i=1}^{N} \sum_{j \neq i}
\E\left(f(\theta_{ij}(\xi))-f(\xi)\right),
\end{equation}
for $f \in C^2_c(\R^N)$, the set of twice continuously differentiable functions with compact support, where $\L_d$ acts on each coordinate, $\L_d f(\xi)=\sum_{i=1}^N\L_d^{(i)}f(\xi)$. Here $\theta_{ij}(\xi)$ is a random transformation of the positions of particle $i$ and $j$, given their respective positions.

Each ordered pair $(i,j)$ thus interacts at rate $1/(N-1)$, so that each
particle initiates interactions at total rate $1$; equivalently, the
tagged-particle clock used throughout the paper is a Poisson process of rate
$1$. See Remark~\ref{rem:normalisation} for the effect of changing this
normalisation.

We denote by $\xi^N(t)=(\xi^N_1(t), \ldots, \xi^N_N(t))$ the process at time $t$ and define the  \emph{empirical measure} associated to $\xi^N$ by
\begin{equation*}
m^N_t:= \frac{1}{N} \sum_{i=1}^{N} \delta_{\xi_t(i)}, \qquad t>0.
\end{equation*}
The \emph{tail of the cumulative distribution function} of $m^N_t$ is given by,
\begin{equation*}
\psi_N(x,t):= m^N_t([x,\infty))	, \qquad x\in \R,t>0.
\end{equation*}
For simplicity, we omit the dependence of $m^N_t$ and $\psi_N(x,t)$ on ${\xi^N}$ if not necessary. In the last decade, many papers have proven the convergence toward a limiting equation of the form:

\be \label{eq:HLgen}
\frac{\partial \psi(x,t)}{\partial t} = \L_d^* \psi(x,t) +  F(\psi)(x,t),
\ee
where $F$ corresponds to a limiting non-linear transformation accounting for the reaction term and $\L_d^*$ is the dual of $\L_d$ (see \cite{bertoin1996}).

\

We now provide a number of examples of such mechanisms and pointers towards  convergence results in the literature.

\begin{example}[Fisher-Kolmogorov-Petrovskii-Piskunov]\label{ex:FKPP}
\begin{equation}
\label{eq:theta_FKPP}
\theta_{ij}(\xi)(k)=
\begin{cases}
\xi(j) & \text{ if $k=i$ and  $\xi(i) < \xi(j)$},\\
\xi(k) & \text{ otherwise},
\end{cases}
\end{equation}
Namely, given two labels $i,j$, $\theta_{ij}$ pushes the position of the leftmost of the two particles to the right-most.
When $\L_d$ is the usual Laplacian operator, \cite{GJM} proved that the limiting system corresponds to the F-KPP equation
\begin{equation}
 \label{KPP}
\begin{array}{l}
\displaystyle {\partial_t u}= \frac{1}{2}  \displaystyle {\partial^2_{x} u} + u^2 -u, \quad x \in \R,\,
t>0,\\
\\
\displaystyle u(0,x)= \displaystyle u_0(x), \quad x \in \R,
\end{array}
\end{equation}
where $u = 1 - \psi$ is the cumulative distribution function associated to the survival function~$\psi$.
Observe that the reaction term of \eqref{KPP} corresponds precisely to the
normalisation of \eqref{eq:generator}: a particle at $y$ overtakes the level
$x>y$ at rate $\psi(x)$, so that $\partial_t\psi=\frac{\sigma^2}{2}\partial^2_x\psi+\psi(1-\psi)$.
Fisher and Kolmogorov, Petrovskii and Piskunov \cite{Fisher,KPP} proved
independently that \eqref{KPP} admits infinitely many travelling waves, one
for each velocity above a critical one; their domains of attraction were
characterised in \cite{AW,U}, and the PDE literature is surveyed in
\cite{KPP,AW}. On the probabilistic side, McKean \cite{MK} exhibited the exact
connection with branching Brownian motion, and Bramson \cite{Bramson78,B2}
established convergence to the wave. A large subsequent literature builds
microscopic models for the front: \cite{BD2,BD} on the influence of
microscopic fluctuations, \cite{BDMM,MMQ,MMQ2,BenDep,DP} adding noise or
cut-off mechanisms, and \cite{BDMM2,DR,M,deMasi,GJM} constructing interacting
particle systems with a conserved number of particles --- the class considered
here.

\end{example}

\begin{example}[Model with independent jumps and copying (Baryshnikov-Stolyar model)]\label{ex:mixed}

The model considered in \cite{baryshnikov2025large} consists of $N$ particles on the real line with two competing mechanisms: at rate $\lambda > 0$, each particle performs an independent jump with size drawn from a distribution with finite mean (which we assume w.l.o.g. to be 1); at rate $\alpha > 0$, particles interact according to the same mechanism as in Example~\ref{ex:FKPP}:
\begin{equation*}
\theta_{ij}(\xi)(k) =
\begin{cases}
\xi(j) & \text{if } k = i \text{ and } \xi(i) < \xi(j), \\
\xi(k) & \text{otherwise},
\end{cases}
\end{equation*}

Throughout we take $\alpha=1$, in accordance with the normalisation of
\eqref{eq:generator}; the general case reduces to this one by rescaling time,
which replaces $\lambda$ by $\lambda/\alpha$ and the speed by $c/\alpha$.

The convergence of the empirical measure to the mean-field limit has been established in \cite{baryshnikov2025large}. The paper proves that, for any fixed $N$, the system has an appropriately defined steady state and the steady-state speed $v_N$ converges to a critical value $v^{*}$ determined by a minimum speed selection principle. The authors also show that there are traveling wave solutions for speeds $v \geq v^{*}$ in the case of exponentially distributed jumps. 

\end{example}

Note that in the example above, one can think of all particles moving according to a compound Poisson process between moments of interaction. For Lévy processes other than Brownian motion and compound Poisson, existence and uniqueness of travelling waves have been studied in \cite{Ky} and \cite{GJ2018}.

\begin{example}[Quantile mechanisms and power of 2]\label{ex:power2}

The model considered in this line of work consists of \( N \) particles on the real line, experiencing "jump urges"  at exponential times. The key feature is the mean-field interaction: the probability that a particle actually jumps upon an urge is a decreasing function \( \eta(v) \) of its quantile \( v \) within the current empirical distribution of all particle locations. When a jump occurs, its size is given by an i.i.d. non-negative random variable \( X \).

A specific example of those mechanisms is the so-called power of 2, which can be described as follows. There is a Poisson process associated to each particle. At the instances of the process, it chooses uniformly between all other particles and jumps (according to a fixed distribution) only if the position of the
 chosen particle is bigger than its own position. In other words:

\begin{equation}
\label{eq:theta_power2}
\theta_{ij}(\xi)(k)=
\begin{cases}
\xi(i) + X & \text{ if $k=i$ and  $\xi(i) \le \xi(j)$},\\
\xi(k) & \text{ otherwise},
\end{cases}
\end{equation}
with $X$ a random variable with a finite first moment. We will assume, w.l.o.g., that $\E(X)=1$. Namely, given two labels $i,j$, $\theta_{ij}$ pushes the position of the leftmost particle by $X$.

The convergence of the empirical measure was established in
\cite{Greenberg1995} under a finite-dependence assumption, then in
\cite{Balazs2014} for a related model and in \cite{Stolyar2023} in general. As
for the travelling wave, \cite{Greenberg1996} showed that if one exists it is
unique up to a shift and attracts the mean-field trajectories; explicit forms
were obtained in special cases in
\cite{Balazs2014,Hongler2015,Hongler2019}; and \cite{Stolyar2023} settled
existence under $\mathbb EX^2<\infty$, with uniqueness under an additional
absolute-continuity assumption on the jump law. None of these results
describes the profile itself, which is what Section~\ref{sec:power-of-2}
provides.

\end{example}

\subsection{Travelling waves as fixed points}\label{sec:tw-fp}

A travelling wave of \eqref{eq:HLgen} is a solution of the form
$\psi(x,t)=h(x-ct)$, $h(-\infty)=1$, $h(\infty)=0$, which transforms
\eqref{eq:HLgen} into
\begin{equation}\label{eq:TWpde}
-ch'(z)=\mathcal L_d h(z)+F(h)(z).
\end{equation}
The point of view we adopt here is to characterise the travelling wave,
viewed as a probability distribution, by a distributional equation. To this
end we write the mean-field equation in weak form, as an equation for the
one-particle law $\rho_t$; \eqref{eq:HLgen} is the form it takes for the
survival function $\psi(x,t)=\rho_t([x,\infty))$.

We assume throughout that the random transformation $\theta_{ij}$ modifies
the position of particle $i$ only, and that it is translation covariant,
\begin{equation*}
\theta_{ij}(\xi+a\mathbf 1)=\theta_{ij}(\xi)+a\mathbf 1,
\qquad a\in\mathbb R;
\end{equation*}
with the notation $\bar\theta(\xi_1,\xi_J):=\theta_{1J}(\xi)(1)-\xi_1$, the
latter property reads $\bar\theta(x+a,y+a)\overset{d}{=}\bar\theta(x,y)$. Both
properties hold in Examples~\ref{ex:FKPP}--\ref{ex:power2}. Formally, the law
of $\bar\theta(x,y)$ is given by a measurable probability kernel. We let
$\alpha>0$ denote the rate at which each particle initiates interactions
($\alpha=1$ under the normalisation of \eqref{eq:generator}, see
Remark~\ref{rem:normalisation}), and, for a probability measure $\nu$ on
$\mathbb R$, we define the Markov operator
\[
K_\nu g(x):=\E\bigl[g\bigl(x+\bar\theta(x,\chi_2)\bigr)\bigr],
\qquad\chi_2\sim\nu,
\]
acting on bounded measurable functions $g$.

We write $C_0(\mathbb R)$ for the continuous functions vanishing at
infinity, with the supremum norm, $C^2_0$ for the twice continuously
differentiable functions $f$ with $f,f',f''\in C_0(\mathbb R)$, and
$C_c^\infty$ for the smooth functions with compact support. For a L\'evy
process $Y$, $\mathcal A_Y$ denotes the generator of its Feller semigroup on
$C_0(\mathbb R)$ and $D(\mathcal A_Y)$ its domain; by
\cite[Thm.~31.5]{Sato99}, $C^2_0\subset D(\mathcal A_Y)$ and $\mathcal A_Y$
maps $C^2_0$ into $C_0(\mathbb R)$. The same holds for $\mathcal L_d$.

\begin{definition}[Mean-field equation, travelling waves]\label{def:TW}
A family $(\rho_t)_{t\ge0}$ of probability measures on $\mathbb R$ is a weak
solution of the mean-field equation if, for every
$g\in C^2_0$ and $t\ge0$,
\begin{equation}\label{eq:weakMF}
\langle\rho_t,g\rangle-\langle\rho_0,g\rangle
=\int_0^t\bigl\langle\rho_s,\ \mathcal L_d g+\alpha\,(K_{\rho_s}g-g)\bigr\rangle\,ds .
\end{equation}
A probability measure $\nu$ is a \emph{travelling wave with speed $c$} if
$\rho_t:=\nu(\cdot-ct)$, i.e.\ the law of $\chi+ct$ with $\chi\sim\nu$, is a
weak solution of \eqref{eq:weakMF}.
\end{definition}


We can now make explicit the link between travelling waves and distributional equations.

\begin{proposition}[Travelling waves as fixed points]\label{prop:TW-FP}
Let $c\in\mathbb R$, let $Y(t):=\Gamma(t)-ct$, with generator
$\mathcal A_Y=\mathcal L_d-c\,\partial_x$, and let $T\sim\mathrm{Exp}(\alpha)$ be
independent of $Y$. For a probability measure $\nu$ on $\mathbb R$ the
following are equivalent:
\begin{enumerate}
\item[(a)] $\nu$ is a travelling wave with speed $c$;
\item[(b)] \textup{(stationarity in the wave frame)} for every
$g\in C^2_0$,
\[
\bigl\langle\nu,\ \mathcal A_Yg+\alpha\,(K_\nu g-g)\bigr\rangle=0;
\]
\item[(c)] \textup{(fixed point)}
\begin{equation}\label{eq:tw-dist}
\chi\;\overset{d}{=}\;\chi_1+\bar\theta(\chi_1,\chi_2)+Y(T),
\end{equation}
where $\chi,\chi_1,\chi_2\sim\nu$, and $\chi_1$, $\chi_2$, the randomness of
$\bar\theta$, $T$ and $Y$ are mutually independent.
\end{enumerate}
If moreover $\E|\chi|<\infty$, $\E|\Gamma(1)|<\infty$ and
$\E|\bar\theta(\chi_1,\chi_2)|<\infty$, then
\[
c=\E\,\Gamma(1)+\alpha\,\E\,\bar\theta(\chi_1,\chi_2).
\]
\end{proposition}

\begin{proof}
The time change $t\mapsto t/\alpha$ turns \eqref{eq:weakMF}, with generator
$\mathcal L_d$, rate $\alpha$ and speed $c$, into the same equation with
generator $\alpha^{-1}\mathcal L_d$, rate $1$ and speed $c/\alpha$. It
replaces $\mathcal A_Y$ by $\alpha^{-1}\mathcal A_Y$, which divides the
left-hand side of (b) by $\alpha$, and it leaves the law of $Y(T)$ unchanged,
since $Y(T)=\tilde Y(\alpha T)$ with $\tilde Y(t):=Y(t/\alpha)$ and
$\alpha T\sim\mathrm{Exp}(1)$. For the equivalence of (a), (b) and (c) we may
therefore assume $\alpha=1$. For $g\in C^2_0$ put
$\Phi(g):=\langle\nu,\mathcal A_Yg+K_\nu g-g\rangle$ and
$g_s:=g(\cdot+cs)$.

\smallskip
Let us first prove \emph{(a)$\Leftrightarrow$(b).} Let $\rho_s=\nu(\cdot-cs)$ and
$\chi\sim\nu$. Since $\mathcal L_d$ commutes with translations,
$\langle\rho_s,\mathcal L_dg\rangle=\langle\nu,\mathcal L_dg_s\rangle$, and,
by translation invariance of $\bar\theta$,
\[
\langle\rho_s,K_{\rho_s}g\rangle
=\E\bigl[g\bigl(\chi+cs+\bar\theta(\chi+cs,\chi_2+cs)\bigr)\bigr]
=\langle\nu,K_\nu g_s\rangle .
\]
Moreover
$\langle\rho_t,g\rangle-\langle\rho_0,g\rangle=\E[g(\chi+ct)-g(\chi)]
=\int_0^tc\,\langle\nu,g_s'\rangle\,ds$. Hence $(\rho_t)$ solves
\eqref{eq:weakMF} if and only if $\int_0^t\Phi(g_s)\,ds=0$ for all $t\ge0$ and
$g\in C^2_0$. The map $s\mapsto\Phi(g_s)$ is continuous: $\mathcal A_Yg$
belongs to $C_0(\mathbb R)$ and is therefore uniformly continuous, and
$K_\nu g_s\to K_\nu g_{s_0}$ boundedly pointwise as $s\to s_0$. Thus (a)
holds if and only if $\Phi(g_s)=0$ for all $s$ and $g$, which is (b) since
$C^2_0$ is invariant under translations.

\smallskip
\emph{Next, (b)$\Leftrightarrow$(c).} Let
$Ru(x):=\E\,u(x+Y(T))=\int_0^\infty e^{-t}\,\E\,u(x+Y(t))\,dt$ for bounded
measurable $u$. Conditioning on $\chi_1+\bar\theta(\chi_1,\chi_2)$, whose law
is $\nu K_\nu$, shows that (c) is equivalent to
\begin{equation}\label{eq:FP-op}
\langle\nu,u\rangle=\langle\nu,K_\nu Ru\rangle
\end{equation}
holding for all $u$ in any class of bounded functions which separates
probability measures.

For $g\in C^2_0$ put $u:=g-\mathcal A_Yg\in C_0(\mathbb R)$. Since $R$ is the
resolvent at $1$ of the Feller semigroup of $Y$, $R=(I-\mathcal A_Y)^{-1}$ on
$C_0(\mathbb R)$, so that $g=Ru$ and
$\mathcal A_Yg+K_\nu g-g=K_\nu Ru-u$: condition (b) for $g$ is exactly
\eqref{eq:FP-op} for $u$. As $g$ ranges over $C^2_0$, $u$ ranges over a class
$\mathcal U\subset C_0(\mathbb R)$ which contains $C^2_0$, and therefore
separates probability measures: for $u\in C^2_0$, the function $g:=Ru$ lies
in $C^2_0$, because $R$ acts by convolution and thus commutes with
differentiation, and $(I-\mathcal A_Y)g=u$. Hence (b) is equivalent to
\eqref{eq:FP-op} on $\mathcal U$, that is, to (c).

\smallskip
\emph{Speed.} Under the moment assumptions, take expectations in
\eqref{eq:tw-dist}, use $\E\,Y(T)=\E\,T\cdot\E\,Y(1)=(\E\,\Gamma(1)-c)/\alpha$,
and cancel $\E\chi$.
\end{proof}

\begin{remark}[The tagged particle]\label{rem:tagged}
Condition (c) has a direct probabilistic meaning. Start the nonlinear
(McKean--Vlasov) process from $\nu$ and follow a tagged particle in the frame
moving at speed $c$; by (b) its law is stationary. Between two consecutive
interactions it initiates, it moves as $Y$ during an $\mathrm{Exp}(\alpha)$
time, and at each of these interactions it is displaced by $\bar\theta$
against a partner sampled from $\nu$; comparing its positions just before two
consecutive interaction epochs gives \eqref{eq:tw-dist}. In what follows,
$\chi$ denotes this tagged particle seen in the wave frame. No assumption on
the underlying particle system is involved: Sections~\ref{sec:prob-FKPP},
\ref{sec:power-of-2} and~\ref{sec:mixed} are statements about travelling waves
of the mean-field equation, obtained by analysing \eqref{eq:tw-dist} for the
corresponding $\bar\theta$.
\end{remark}

\begin{remark}[Identification of the speed]\label{rem:speed}
For $\alpha=1$ and whenever $\E|\chi|$ and the two expectations are finite,
Proposition~\ref{prop:TW-FP} identifies the speed as
\begin{equation}\label{eq:speed-identity}
c=\gamma_1+\bar\beta,
\qquad \gamma_1:=\mathbb E\,\Gamma(1),
\quad \bar\beta:=\mathbb E_{h^{\otimes2}}\,\bar\theta(\chi_1,\chi_2),
\end{equation}
consistently with the fact that the barycenter drifts at the mean rate at
which diffusion and interactions displace mass. For the synchronisation
mechanism this reads
$\bar\beta=\mathbb E\,(\chi_2-\chi_1)^+=\tfrac12\,\mathbb E\,|\chi_1-\chi_2|$.
For the power-of-2 mechanism, $\bar\theta(x,y)=X\,\mathbf 1\{x\le y\}$ with $X$
independent of $(\chi_1,\chi_2)$, so that, whenever $h$ is atomless,
\begin{equation}\label{eq:speed-pow2}
\bar\beta=\mathbb E[X]\;\mathbb P(\chi_1\le \chi_2)=\frac{\mathbb E[X]}{2},
\qquad\text{so that}\qquad c=\gamma_1+\frac{\mathbb E[X]}{2}.
\end{equation}

Note that \eqref{eq:speed-identity} is a compatibility relation between the
speed and the profile rather than a formula determining $c$: for
synchronisation $\bar\beta$ depends on $h$, and \eqref{eq:speed-identity}
holds identically along the whole one-parameter family of travelling waves,
selecting none of them (Example~\ref{ex:pure_copying}). Which member is
selected depends on how the wave is reached: for the mean-field dynamics it
is determined by the initial condition \cite{AW,U}, and for the stationary
finite-$N$ systems it is fixed by hypothesis (ii) of Theorem~\ref{prop:tw}
(see Remark~\ref{rem:selection}). Relation \eqref{eq:speed-identity}
\emph{does} determine $c$ when $\bar\theta(x,y)$ depends on $(x,y)$ only
through the sign of $x-y$, since then $\bar\beta$ takes the same value for
every atomless $h$ by exchangeability. The power-of-2 mechanism is of this
type, which is why \eqref{eq:speed-pow2} is profile-free and why
Section~\ref{sec:power-of-2} may fix the speed at the outset.
\end{remark}

\begin{remark}[Normalisation of the interaction clocks]\label{rem:normalisation}
We take $\alpha=1$, which puts $T\sim\mathrm{Exp}(1)$ in \eqref{eq:tw-dist},
except in Section~\ref{sec:mixed} where $\alpha=\lambda$. This is the
normalisation under which Example~\ref{ex:FKPP} produces exactly \eqref{KPP},
and under which the model of Example~\ref{ex:mixed} is described in
\cite{baryshnikov2025large}. For a general $\alpha$,
\eqref{eq:speed-identity} becomes $c=\gamma_1+\alpha\bar\beta$; the
power-of-2 speed is then $\alpha\,\mathbb E[X]/2$, so that the value
$c=\mathbb E[X]$ of \cite{Stolyar2023} corresponds to $\alpha=2$. A change of
$\alpha$ also changes the ratio between the diffusive and the interaction time
scales, so the two conventions are not related by a time change.
\end{remark}

\subsection{From particle systems to travelling waves}\label{sec:particles}

Proposition~\ref{prop:TW-FP} concerns the mean-field equation alone. We now
connect it to the particle systems of Section~\ref{sec:mf}, through the
stationary regime of the finite-$N$ systems. For finitely many particles the
finite-$N$ system travels at a speed $v_N$ and, even after subtracting the
correct linear drift, the configuration is not anchored to any fixed origin.
We therefore consider the process, stationary under suitable recurrence
assumptions, describing the configuration seen from its centre of mass, and
introduce the barycenter and the recentered process
\begin{equation*}
\bar\xi^N(t):=\frac1N\sum_{i=1}^N\xi_i^N(t),
\qquad
\widehat Z_i^N(t):=\xi_i^N(t)-\bar\xi^N(t),\qquad 1\le i\le N.
\end{equation*}
Translation covariance of $\theta_{ij}$ implies that
$\widehat Z^N=(\widehat Z^N_1,\dots,\widehat Z^N_N)$ is itself a Markov
process, with values in the hyperplane
$\mathcal H_N=\{z\in\mathbb R^N:\sum_{i}z_i=0\}$.

\begin{theorem}
[Travelling waves via a tagged particle]\label{prop:tw}
Assume that:
\begin{enumerate}
\item[(i)] for every $N$, the recentred process $\widehat Z^N$ admits a
stationary distribution $\pi^N$, and the system is in the corresponding
stationary regime;
\item[(ii)] \textup{(propagation of chaos in the stationary regime)} for
every fixed $k$, the law of $(\widehat Z^N_1,\dots,\widehat Z^N_k)$ under
$\pi^N$ converges weakly, as $N\to\infty$, to $h^{\otimes k}$ for some
probability measure $h$ on $\mathbb R$;
\item[(iii)]
\textup{(Law of large numbers for the barycenter)} there exists
$c\in\mathbb R$ such that, in the stationary regime, for every $t\ge0$,
\begin{equation*}
\sup_{s\le t}\bigl|\bar\xi^N(s)-\bar\xi^N(0)-cs\bigr|
\xrightarrow[N\to\infty]{\;\mathbb P\;}0;
\end{equation*}
\item[(iv)] $\bar\theta$ is continuous $h\otimes h$-a.s.
\end{enumerate}
Then $h$ satisfies \eqref{eq:tw-dist} with $\alpha=1$, that is,
\[
\chi\;\overset{d}{=}\;\chi_1+\bar\theta(\chi_1,\chi_2)+\Gamma(T)-cT,
\qquad T\sim\mathrm{Exp}(1),
\]
with the independence structure of Proposition~\ref{prop:TW-FP}(c).
Consequently, $h$ is a travelling wave with speed $c$ of the mean-field
equation \eqref{eq:weakMF}. Furthermore, if $\E|\chi|<\infty$,
$\mathbb E\,|\Gamma(1)|<\infty$ and
$\mathbb E\,|\bar\theta(\chi_1,\chi_2)|<\infty$, then
\begin{equation*}
c=\mathbb E\,\Gamma(1)+\mathbb E\,\bar\theta(\chi_1,\chi_2).
\end{equation*}
\end{theorem}

\begin{proof}
We shall work with the stationary version of the process provided by (i). Since each particle
initiates an interaction with each of the other $N-1$ particles at rate
$1/(N-1)$, the epochs at which particle $1$ initiates an interaction form a
Poisson process of rate $1$, independent of the driving L\'evy noises and of
the interaction clocks of the other particles. Let $\tau_0<\tau_1$ be two
consecutive such epochs, set $T:=\tau_1-\tau_0$, and let $J$ denote the
(uniformly chosen) partner at time~$\tau_0$. By memorylessness,
$T\sim\mathrm{Exp}(1)$ is independent of
$\bigl(\widehat Z^N(\tau_0-),J\bigr)$ and of the L\'evy increment appearing
below.

\smallskip

Since the interaction times of particle $1$ form a Poisson process independent of everything else, the PASTA property yields
\begin{equation}\label{eq:pasta}
\widehat Z^N(\tau_0-)\sim\pi^N
\qquad\text{and}\qquad
\widehat Z^N(\tau_1-)\sim\pi^N.
\end{equation}
It is important here that we compare the states \emph{just before} the
interaction epochs: the time-stationary law is the one seen by the Poisson
clock, whereas the post-jump law need not coincide with $\pi^N$.

\smallskip
Between $\tau_0$ and $\tau_1$ particle $1$ initiates no interaction and, since
$\theta_{ij}$ displaces the initiating particle only, its position evolves
only through its own L\'evy noise. Hence, writing
$\bar\theta_0:=\bar\theta\bigl(\xi^N_1(\tau_0-),\xi^N_J(\tau_0-)\bigr)$,
\begin{equation}\label{eq:onestep}
\xi^N_1(\tau_1-)=\xi^N_1(\tau_0-)+\bar\theta_0+\Gamma(T),
\end{equation}
where, conditionally on $T$, $\Gamma(T)$ is the L\'evy displacement over an
interval of length $T$, independent of everything else. By translation
invariance of $\bar\theta$ we may recentre its arguments:
$\bar\theta_0=\bar\theta\bigl(\widehat Z^N_1(\tau_0-),
\widehat Z^N_J(\tau_0-)\bigr)$.

The interaction at time $\tau_0$ displaces the barycenter by
$\bar\theta_0/N$, that is,
$\bar\xi^N(\tau_0)=\bar\xi^N(\tau_0-)+\bar\theta_0/N$. Subtracting
$\bar\xi^N(\tau_1-)$ from both sides of \eqref{eq:onestep} and adding and
subtracting $\bar\xi^N(\tau_0-)$ therefore gives
\begin{equation}\label{eq:recentred}
\widehat Z^N_1(\tau_1-)
=\widehat Z^N_1(\tau_0-)
+\Bigl(1-\frac1N\Bigr)\bar\theta_0
+\Gamma(T)
-\bigl(\bar\xi^N(\tau_1-)-\bar\xi^N(\tau_0)\bigr).
\end{equation}

\smallskip
We can now consider the passage to the limit.
By \eqref{eq:pasta}, both $\widehat Z^N_1(\tau_1-)$ and
$\widehat Z^N_1(\tau_0-)$ are distributed according to the first marginal of
$\pi^N$. By (ii), exchangeability and the uniform choice of $J$,
\begin{equation*}
\bigl(\widehat Z^N_1(\tau_0-),\widehat Z^N_J(\tau_0-)\bigr)
\;\Longrightarrow\;(\chi_1,\chi_2)\sim h\otimes h.
\end{equation*}
By (iii), applied from the time origin $\tau_0$ (see
Remark~\ref{rem:palm} below), and by conditioning on the independent
duration $T$ and dominated convergence,
\begin{equation*}
\bar\xi^N(\tau_1-)-\bar\xi^N(\tau_0)-cT
\xrightarrow[N\to\infty]{\;\mathbb P\;}0 .
\end{equation*}
Finally, $(1-1/N)\,\bar\theta_0\Rightarrow\bar\theta(\chi_1,\chi_2)$ by (iv) and the
continuous mapping theorem. Applying Slutsky's theorem to
\eqref{eq:recentred} yields \eqref{eq:tw-dist} with $\alpha=1$. The last two
assertions then follow from Proposition~\ref{prop:TW-FP}.
\end{proof}

\begin{remark}[Shift of the time origin in (iii)]\label{rem:palm}
Assumption (iii) is formulated in the time-stationary regime, whereas in the
passage to the limit above, the barycenter increment is taken from the epoch
$\tau_0$, which is
distributed according to the Palm measure of particle $1$'s clock. The
discrepancy between the two only enters through the trajectory of particle
$1$ itself, whose contribution to the barycenter increment over the window
$(\tau_0,\tau_1)$ is $O_{\mathbb P}(1/N)$, and through the particles that
interact with particle $1$ during this window: each of the other $N-1$
particles chooses particle $1$ at rate $1/(N-1)$, so only $O_{\mathbb P}(1)$
of them are affected, again contributing $O_{\mathbb P}(1/N)$ to the average.
Alternatively, we might simply have stated (iii) with the time origin at $\tau_0$,
i.e.\ under the Palm measure.
\end{remark}

\begin{remark}[Selection]\label{rem:selection}
Proposition~\ref{prop:TW-FP} and Theorem~\ref{prop:tw} play different roles.
The former characterises \emph{all} travelling waves of the mean-field
equation, at every admissible speed. The latter concerns the order of limits
$t\to\infty$ then $N\to\infty$: under (i)--(iv) the stationary finite-$N$
systems single out one wave $h$, whatever their initial condition. Where a
selection principle is known --- \cite{GJM} for Brownian synchronisation, and
\cite{baryshnikov2025large} for the speeds in the compound Poisson case ---
this is the wave of minimal speed. The other waves are reached in the reverse
order of limits, by the mean-field dynamics started from suitable initial
conditions \cite{AW,U}.
\end{remark}


We now show that, under the standing assumptions of
Theorem~\ref{prop:tw}, the law of large numbers for the barycenter
(assumption~(iii)) is not independent of the others: it follows from
(i)--(ii) together with second-moment bounds, and the speed $c$ is
identified intrinsically in terms of the limiting profile $h$. Recall that
$\bar\theta(x,y)$ denotes the (random) displacement of the initiating
particle; formally, its law is given by a measurable kernel, and we set
\begin{equation*}
\beta(x,y):=\mathbb E\bigl[\bar\theta(x,y)\bigr],
\qquad
\beta_2(x,y):=\mathbb E\bigl[\bar\theta(x,y)^2\bigr].
\end{equation*}
Note that propagation of chaos alone cannot control the barycenter: its
drift is carried by an unbounded functional of the empirical measure, and
weak convergence does not see it (a configuration of vanishing probability
with one particle at distance of order $N$ leaves chaos intact but shifts
the barycenter drift by a quantity of order one). Uniform integrability
enters through the following two assumptions:
\begin{align}
&\tag{M}\label{eq:M}
C:=\sup_{N}\ \mathbb E_{\pi^N}\bigl[(\widehat Z^N_1)^2\bigr]<\infty ;\\
&\tag{G}\label{eq:G}
\beta_2(x,y)\le C_\beta\bigl(1+x^2+y^2\bigr)
\quad\text{for some }C_\beta<\infty\text{ and all }x,y\in\mathbb R .
\end{align}
We write $\widehat m^N_t:=\frac1N\sum_{i=1}^N\delta_{\widehat Z^N_i(t)}$ for
the recentred empirical measure.

\begin{lemma}[LLN for the barycenter, general mechanism]\label{lem:bary-gen}
Consider the particle system with generator \eqref{eq:generator}, where the
transformation $\theta_{ij}$ displaces the initiating particle $i$ only and
is translation covariant, and the driving L\'evy process satisfies
$\gamma_1:=\mathbb E\,\Gamma(1)\in\mathbb R$ and
$\gamma_2:=\operatorname{Var}\Gamma(1)<\infty$. Assume:
\begin{enumerate}
\item[(i)] for every $N$ the recentred process $\widehat Z^N$ is in an
exchangeable stationary regime $\pi^N$;
\item[(ii)] \textup{(two-particle chaos)} the law of
$(\widehat Z^N_1,\widehat Z^N_2)$ under $\pi^N$ converges weakly to
$h\otimes h$;
\item[(iii)] the moment bounds \eqref{eq:M} and \eqref{eq:G} hold;
\item[(iv)] $\beta$ is continuous $h\otimes h$-almost everywhere.
\end{enumerate}
Then $\mathbb E_h[\chi^2]\le C$, the \emph{mean interaction displacement}
\begin{equation}\label{eq:betabar}
\bar\beta:=\mathbb E_{h^{\otimes2}}\bigl[\bar\theta(\chi_1,\chi_2)\bigr]
 =\mathbb E_{h^{\otimes2}}\bigl[\beta(\chi_1,\chi_2)\bigr]
\end{equation}
and the constant
\begin{equation}\label{eq:cdef-gen}
c:=\gamma_1+\bar\beta
\end{equation}
are finite, and for every $t\ge0$,
\begin{equation*}
\sup_{s\le t}\bigl|\bar\xi^N(s)-\bar\xi^N(0)-cs\bigr|
\xrightarrow[N\to\infty]{\;\mathbb P\;}0 .
\end{equation*}
In particular, assumption (iii) of Theorem~\ref{prop:tw} holds with $c$
given by \eqref{eq:cdef-gen}.
\end{lemma}

\begin{proof}
Between interactions each $\xi_i$ moves as an independent copy $\Gamma_i$ of
the driving L\'evy process, and an $(i,j)$-interaction, occurring at rate
$1/(N-1)$, displaces the barycenter by $\bar\theta(\xi_i,\xi_j)/N$. Since the
noises are independent of the interaction clocks and marks,
\begin{equation}\label{eq:decomp-gen}
\bar\xi^N(t)-\bar\xi^N(0)
=\gamma_1t+G^N(t)+\int_0^t b^N(s)\,ds+M^N(t),
\end{equation}
where $G^N(t):=\frac1N\sum_{i=1}^N\bigl(\Gamma_i(t)-\gamma_1t\bigr)$ is a
square-integrable martingale with variance $\gamma_2t/N$; the drift is, by
translation invariance of $\bar\theta$,
\begin{equation}\label{eq:drift-gen}
b^N(s)
=\frac{1}{N(N-1)}\sum_{i\neq j}
\beta\bigl(\widehat Z^N_i(s),\widehat Z^N_j(s)\bigr);
\end{equation}
and $M^N$ is the compensated sum of the interaction jumps, with predictable
quadratic variation
\begin{equation*}
\langle M^N\rangle_t
=\int_0^t\frac{1}{N^2(N-1)}\sum_{i\neq j}
\beta_2\bigl(\widehat Z^N_i(s),\widehat Z^N_j(s)\bigr)\,ds .
\end{equation*}
By stationarity, exchangeability, \eqref{eq:G} and \eqref{eq:M},
\begin{equation}\label{eq:qv-gen}
\mathbb E\,\langle M^N\rangle_t
=\frac{t}{N}\,\mathbb E_{\pi^N}\bigl[\beta_2(\widehat Z^N_1,\widehat Z^N_2)\bigr]
\le\frac{C_\beta(1+2C)\,t}{N};
\end{equation}
in particular the compensator in \eqref{eq:decomp-gen} is integrable and
$M^N$ is a well-defined square-integrable martingale.

\smallskip
\emph{Diffusive and martingale parts.}
By Doob's $L^2$ inequality,
$\mathbb P\bigl(\sup_{s\le t}|G^N(s)|>\varepsilon\bigr)
\le4\gamma_2t/(N\varepsilon^2)\to0$, and, by Doob's $L^2$ inequality
together with \eqref{eq:qv-gen},
\begin{equation*}
\mathbb P\Bigl(\sup_{s\le t}\bigl|M^N(s)\bigr|>\varepsilon\Bigr)
\le\frac{4\,\mathbb E\langle M^N\rangle_t}{\varepsilon^2}
\le\frac{4C_\beta(1+2C)\,t}{N\varepsilon^2}
\xrightarrow[N\to\infty]{}0 .
\end{equation*}

\smallskip
\emph{Finiteness of $c$ and reduction to convergence in probability.}
By (ii), the Skorokhod representation theorem and Fatou's lemma,
$\mathbb E_h[\chi^2]\le\liminf_N\mathbb E_{\pi^N}[(\widehat Z^N_1)^2]\le C$;
hence, by Jensen's inequality and \eqref{eq:G},
$\mathbb E_{h^{\otimes2}}|\beta|
\le(\mathbb E_{h^{\otimes2}}\beta_2)^{1/2}\le\sqrt{C_\beta(1+2C)}<\infty$,
so $\bar\beta$, and hence $c$, is finite. Moreover $b^N(0)$ is an average of
the $N(N-1)$ variables $\beta(\widehat Z^N_i,\widehat Z^N_j)$, $i\neq j$, so
Jensen's inequality, exchangeability, $\beta^2\le\beta_2$, \eqref{eq:G} and
\eqref{eq:M} give
\begin{equation}\label{eq:L2-drift}
\mathbb E\bigl[b^N(0)^2\bigr]
\le\mathbb E_{\pi^N}\bigl[\beta(\widehat Z^N_1,\widehat Z^N_2)^2\bigr]
\le\mathbb E_{\pi^N}\bigl[\beta_2(\widehat Z^N_1,\widehat Z^N_2)\bigr]
\le C_\beta(1+2C).
\end{equation}
The family $\{b^N(0)\}_N$ is therefore bounded in $L^2$, hence uniformly
integrable, and it suffices to prove that
$b^N(0)\to\bar\beta$ \emph{in probability}: uniform integrability then
upgrades this to $\mathbb E|b^N(0)-\bar\beta|\to0$. Note that it is $\bar\beta$
and not $c$ that is the limit here: the diffusive contribution $\gamma_1t$ was
split off in \eqref{eq:decomp-gen} and does not enter the drift $b^N$.

\smallskip
\emph{Drift part.}
Fix $K>0$, set $\beta_K:=(\beta\vee(-K))\wedge K$, let
\begin{equation*}
b^N_K(0):=\frac{1}{N(N-1)}\sum_{i\neq j}
\beta_K\bigl(\widehat Z^N_i,\widehat Z^N_j\bigr),
\qquad
I^N_K:=\iint\beta_K\,d\widehat m^N_0\otimes d\widehat m^N_0,
\end{equation*}
and $\bar\beta_K:=\mathbb E_{h^{\otimes2}}[\beta_K(\chi_1,\chi_2)]$. We compare the
four quantities $b^N(0)$, $b^N_K(0)$, $I^N_K$, $\bar\beta_K$ in turn.

Since $|\beta-\beta_K|=(|\beta|-K)^+\le\beta^2/K\le\beta_2/K$, exchangeability
and \eqref{eq:G}--\eqref{eq:M} give
\begin{equation}\label{eq:trunc1}
\mathbb E\bigl|b^N(0)-b^N_K(0)\bigr|
\le\frac{1}{K}\,\mathbb E_{\pi^N}\bigl[\beta_2(\widehat Z^N_1,\widehat Z^N_2)\bigr]
\le\frac{C_\beta(1+2C)}{K},
\end{equation}
and, using $\mathbb E_h[\chi^2]\le C$,
\begin{equation}\label{eq:trunc2}
\bigl|\bar\beta_K-\bar\beta\bigr|
\le\frac{1}{K}\,\mathbb E_{h^{\otimes2}}\bigl[\beta_2(\chi_1,\chi_2)\bigr]
\le\frac{C_\beta(1+2C)}{K}.
\end{equation}
Next, $b^N_K(0)$ and $I^N_K=\frac1{N^2}\sum_{i,j}\beta_K(\widehat Z^N_i,\widehat Z^N_j)$
differ only through the normalisation of the off-diagonal terms and through
the $N$ diagonal ones; as $|\beta_K|\le K$, this discrepancy is
\emph{deterministically} bounded:
\begin{equation}\label{eq:diag}
\bigl|b^N_K(0)-I^N_K\bigr|
\le\Bigl(\frac{1}{N(N-1)}-\frac{1}{N^2}\Bigr)N(N-1)K+\frac{N K}{N^2}
=\frac{2K}{N}.
\end{equation}
Finally, since $\pi^N$ is exchangeable, two-particle chaos (ii) is
equivalent to $\widehat m^N_0\to h$ weakly in probability
\cite[Proposition~2.2]{Sznitman91}, whence
$\widehat m^N_0\otimes\widehat m^N_0\to h\otimes h$ weakly in probability.
The discontinuity set of $\beta_K$ is contained in that of $\beta$, which is
$h\otimes h$-null by (iv); as $\beta_K$ is bounded, the mapping theorem
(applied along almost surely convergent subsequences) yields
\begin{equation}\label{eq:mapping}
I^N_K\;\xrightarrow[N\to\infty]{\;\mathbb P\;}\;\bar\beta_K .
\end{equation}

Let now $\varepsilon>0$ and choose $K>C_\beta(1+2C)/\varepsilon$, so that the
right-hand side of \eqref{eq:trunc2} is at most $\varepsilon$. Combining
\eqref{eq:trunc1}--\eqref{eq:mapping} with Markov's inequality,
\begin{equation*}
\limsup_{N\to\infty}\;
\mathbb P\bigl(|b^N(0)-\bar\beta|>3\varepsilon\bigr)
\le\frac{C_\beta(1+2C)}{K\varepsilon},
\end{equation*}
and letting $K\to\infty$ gives $b^N(0)\to\bar\beta$ in probability, hence
in $L^1$ by \eqref{eq:L2-drift}.

\smallskip
Finally, $b^N(s)$ is a fixed function of $\widehat Z^N(s)$, whose law does
not depend on $s$ in the stationary regime, so
\begin{equation*}
\mathbb E\Bigl[\sup_{s\le t}
\Bigl|\int_0^s\bigl(b^N(u)-\bar\beta\bigr)\,du\Bigr|\Bigr]
\le t\,\mathbb E\bigl|b^N(0)-\bar\beta\bigr|
\xrightarrow[N\to\infty]{}0 ;
\end{equation*}
no time-ergodicity of the stationary regime is required. Summing the four
parts of \eqref{eq:decomp-gen} concludes the proof.
\end{proof}


\begin{remark}[Verification in the examples]\label{rem:examples-lln}
Conditions \eqref{eq:G} and (iv) are immediate in
Examples~\ref{ex:FKPP}--\ref{ex:power2}. For synchronisation $\beta(x,y)=(y-x)^+$ is continuous with
$\beta_2(x,y)\le2x^2+2y^2$, and $\bar\beta$ takes the value recorded in
Remark~\ref{rem:speed}, which also notes that the lemma identifies $c$ only in
terms of the profile $h$ furnished by (ii); for the Baryshnikov--Stolyar model
the finite variance of the driver amounts to $\mathbb E[X^2]<\infty$.

For power-of-2, $\beta(x,y)=\mathbb E[X]\mathbf1\{x\le y\}$ and
$\beta_2\le\mathbb E[X^2]$, so \eqref{eq:G} holds under the standing
assumption $\mathbb E[X^2]<\infty$ of \cite{Stolyar2023}, and (iv) holds
whenever $h$ is atomless. Here the drift is deterministic and no propagation
of chaos is needed: when the positions are a.s.\ pairwise distinct,
$b^N(s)=\mathbb E[X]/2$ almost surely, while
$\mathbb E\langle M^N\rangle_t=\mathbb E[X^2]t/(2N)\to0$, so that
\eqref{eq:cdef-gen} gives
\begin{equation}\label{eq:pow2-speed}
\bar\beta=\frac{\mathbb E[X]}{2},
\qquad c=\gamma_1+\frac{\mathbb E[X]}{2},
\end{equation}
which is $1/2$ for a centred driver with $\mathbb E[X]=1$, in agreement with
\eqref{eq:speed-pow2}.
\end{remark}

\begin{remark}[On the continuity assumption in Theorem~\ref{prop:tw}]
\label{rem:ae-cont}
Assumption (iv) above is the natural weakening of the continuity of
$\bar\theta$ required in Theorem~\ref{prop:tw}: there too the continuous
mapping theorem only needs $\bar\theta$, in an appropriate realisation, to be
continuous at $h\otimes h$-almost every point. This is what the power-of-2
mechanism requires, for the reason given in Remark~\ref{rem:examples-lln}.
\end{remark}

\section{Distributional equations} \label{sec:distr_equations}

In this Section we consider the non-linear equation of the introduction, written with the
notation used in the rest of the paper,
\begin{equation} \label{eq:main_deq_general}
\chi \overset{d}{=} \max\{\chi_1,\chi_2\} - A,
\end{equation}
where $\chi_1,\chi_2$ are independent copies of $\chi$ and $A$ is independent of
them. The following lemma is the correspondence with the linear equation
announced there; versions of it appear in the literature, for instance in
\cite{aldous2005}, and we include a proof for completeness.

\begin{lemma} \label{lemma:deq_main}
There exists a solution to \eqref{eq:main_deq_general} of the form
$$
F(x) := \Prob(\chi \le x) = \E\left(e^{-Ve^{-\gamma x}}\right)
$$
if and only if there exists a positive random variable $V > 0$ almost surely and a constant $\gamma > 0$ satisfying:
\begin{equation} \label{eq:deq_main_aux}
V \overset{d}{=} (V_1 + V_2) e^{-\gamma A},
\end{equation}
where $V_1$ and $V_2$ are independent copies of $V$, and all random variables on the right-hand side are mutually independent.
\end{lemma}

\begin{proof}
Suppose first that $V>0$ almost surely and $\gamma>0$ satisfy
\eqref{eq:deq_main_aux}, and set $F(x):=\E[e^{-Ve^{-\gamma x}}]$, which is a
distribution function since it increases from $0$ to $1$. Let $\chi_1,\chi_2$ be
independent with law $F$ and let $A$ be independent of them. Conditionally on
$A$, independence of $V_1,V_2$ gives
$F(x+A)^2=\E\bigl[e^{-(V_1+V_2)e^{-\gamma(x+A)}}\bigm|A\bigr]$, whence
\[
\Prob\bigl(\max\{\chi_1,\chi_2\}-A\le x\bigr)=\E\bigl[F(x+A)^2\bigr]
=\E\bigl[e^{-(V_1+V_2)e^{-\gamma A}e^{-\gamma x}}\bigr]
=\E\bigl[e^{-Ve^{-\gamma x}}\bigr]=F(x),
\]
the third equality by \eqref{eq:deq_main_aux}. Thus $\chi$ solves
\eqref{eq:main_deq_general}.

Conversely, let $\chi$ solve \eqref{eq:main_deq_general} with
$F(x)=\E[e^{-Ve^{-\gamma x}}]$ for some $V>0$ and $\gamma>0$. The computation
above shows that $\E[e^{-sV}]=\E[e^{-s(V_1+V_2)e^{-\gamma A}}]$ at every
$s=e^{-\gamma x}$, hence at every $s>0$; by uniqueness of Laplace transforms
of positive random variables, \eqref{eq:deq_main_aux} holds.
\end{proof}

\begin{remark} \label{rem:deq_main}
The Lemma as stated establishes an equivalence between solutions of \eqref{eq:main_deq_general} of the specific exponential form and solutions of \eqref{eq:deq_main_aux}. For the broader claim that existence of a solution to \eqref{eq:main_deq_general} is equivalent to existence of a solution to \eqref{eq:deq_main_aux}, one would need to address whether all solutions of \eqref{eq:main_deq_general} must be of this exponential form.

In order to demonstrate the connection between the two types of equations further, let us note that, as all the random variables on the RHS of \eqref{eq:general_type_sum} are non-negative, we can introduce the Laplace transform of $\chi$ denoted by $\phi(s)=\E\left(e^{-s\chi}\right)$ which is well-defined for all $s\ge 0$ and such that $\phi(0)=1$, $ \phi(\infty)=0$ and, thanks to \eqref{eq:general_type_sum}, satisfies
\begin{equation} \label{eq:general_type_sum_function}
\phi(s)=\int_{a=0}^\infty \Prob(A \in da) \phi^2(sa).
\end{equation}
We will discuss below conditions for there to exist a unique solution to \eqref{eq:general_type_sum}, or equivalently a unique completely monotone function satisfying \eqref{eq:general_type_sum_function}.

Note also that, with notation $\eta=e^{-\chi}$ and $\tilde{B}=e^{-B}$, the equation \eqref{eq:general_type_max} may be rewritten in the equivalent form as
$$
\eta = \min\{\eta_1,\eta_2\}/\tilde{B},
$$
where now all the random variables on the RHS are non-negative. With the notation $\overline{F}(x) = \Prob(\eta>x)$, one can rewrite the above in the equivalent form as
\begin{equation} \label{eq:general_type_max_function}
\overline{F}(x) = \int_{b=0}^\infty \Prob(\tilde{B} \in db) \overline{F}^2(xb).
\end{equation}
One can now observe that equations \eqref{eq:general_type_sum_function} and \eqref{eq:general_type_max_function} are the same, but, when solving the latter, we are looking for a non-increasing function $\overline{F}$ which does not need to be completely monotone. This wider class is the one treated by Alsmeyer, Biggins and Meiners \cite{ABM12}, who describe the fixed points of $f(t)=\E\prod_i f(tT_i)$ among all non-increasing functions $f$ with $f(0)=1$, not only Laplace transforms. In the non-lattice case, which covers all the variables $A$ arising in this paper, their description should yield exactly the Gumbel mixtures of Lemma~\ref{lemma:deq_main}, hence uniqueness of the wave among all solutions of \eqref{eq:main_deq_general}; we do not check their hypotheses in detail here.
\end{remark}

\begin{remark}
The condition that $V > 0$ almost surely is crucial. If $V$ has atoms at 0, then the representation $F(x) = \E\left(e^{-Ve^{-\gamma x}}\right)$ would give $\lim_{x \to \infty} F(x) = \Prob(V = 0) < 1$, which is not a valid distribution function.
\end{remark}

\subsection{Existence and uniqueness conditions for the linear equation}

Equation~\eqref{eq:deq_main_aux} belongs to a well-studied class of linear distributional equations.
In particular, many papers from the branching literature (among which the seminal papers \cite{Biggins1,Biggins2,Bigginsetal91}) have given conditions of existence and uniqueness of solutions to this equation.

The following proposition summarises the existence and uniqueness results of a solution of \eqref{eq:deq_main_aux} with {\bf finite and with infinite mean}. It corresponds to results of
\cite[Proposition~2.2(ii)]{liu1999asymptotic}
for the case where $ \mathbb{E}\left[A\,e^{-\gamma A}\right] > 0$ and \cite{Liu2001}
 in the critical case when $\mathbb{E}\left[A\,e^{-\gamma A}\right] = 0$.

\begin{proposition}[Existence and uniqueness for the linear equation]\label{prop:linear-existence}
Let \(A\) be a real-valued random variable and \(\gamma>0\) be such that
\(\mathbb{E}\left[e^{-\gamma A}\right]<\infty\) and
\(\mathbb{E}\left[|A|\,e^{-\gamma A}\right]<\infty\).
\begin{itemize}
    \item {\bf Supercritical speed:}
    There exists a unique nonnegative random variable (up to multiplicative constant) \(V\) with \(\mathbb{E}[V]<\infty\) that satisfies \eqref{eq:deq_main_aux} if and only if
\begin{align}
    \mathbb{E}\left[e^{-\gamma A}\right] &= \frac{1}{2}, \label{eq:cond1} \\
    \mathbb{E}\left[A\,e^{-\gamma A}\right] &> 0. \label{eq:cond2}
\end{align}
\item {\bf Critical speed:}
 There exists a unique (up to multiplicative constant) nonnegative random variable \(V\) with \(\mathbb{E}[V]=\infty\) that satisfies \eqref{eq:deq_main_aux} if
there exists $\delta >0$
\begin{align}
    \mathbb{E}\left[e^{-\gamma A}\right] &= \frac{1}{2}, \label{eq:cond1ct} \\
    \mathbb{E}\left[A\,e^{-\gamma A}\right] = 0, \  &  \mathbb{E}\left[|A|^{1+\delta}\,e^{-\gamma A}\right] < \infty. \label{eq:cond2ct}
\end{align}
    \end{itemize}
\end{proposition}

\begin{remark}[On the sign of $A$]\label{rem:sign-of-A}
The variable $A$ is \emph{not} assumed nonnegative, and it is not in the
applications of Section~\ref{sec:prob-FKPP}, where $A=cT-\Gamma(T)$ takes
both signs. In Section~\ref{sec:power-of-2}, by contrast, $A=\overline X+M>0$ almost
surely, which is what makes the supercritical condition \eqref{eq:cond2}
automatic there and places the power-of-2 model always in the supercritical
regime.
\end{remark}

\section{Travelling waves for the synchronisation mechanism}\label{sec:prob-FKPP}

In this section we consider the model where each particle moves according to a L\'evy process $\Gamma$ with cumulant
generating function $\kappa(\theta)=\log\E[e^{\theta\Gamma(1)}]$, and at each
interaction the lower of the two interacting particles jumps to the position
of the higher. Here $\bar\theta(x,y)=(y-x)^+$, so that
$\chi_1+\bar\theta(\chi_1,\chi_2)=\max\{\chi_1,\chi_2\}$ and Proposition~\ref{prop:TW-FP} (with $\alpha=1$) gives
\begin{equation} \label{eq:main_deq_kpp}
\chi \overset{d}{=} \max\{\chi_1,\chi_2\} - A,
\qquad A:=cT-\Gamma(T),\quad T\sim\mathrm{Exp}(1).
\end{equation}
This is precisely \eqref{eq:main_deq_general}, and the whole of the present
section is Proposition~\ref{prop:linear-existence} applied to it. All that is
needed is therefore the transform $\theta\mapsto\E[e^{-\theta A}]$, which we
compute first.

By Proposition~\ref{prop:TW-FP}, the results of this section are statements
about the mean-field equation, valid for an arbitrary L\'evy driver. The
mean-field limit itself has been established in two cases: Brownian motion
(Example~\ref{ex:FKPP}, leading to the F-KPP equation) and the compound
Poisson process (Example~\ref{ex:mixed}, leading to the model of
\cite{baryshnikov2025large}). These are treated in
Sections~\ref{sec:brownian} and~\ref{sec:BS2024}; we conjecture that the
mean-field limit extends to a general L\'evy driver, but do not pursue this
here.

\subsection{Transform calculation}\label{sec:BRW}

Conditioning on $T$ and using $\E[e^{\theta\Gamma(t)}]=e^{t\kappa(\theta)}$,
\begin{equation}\label{eq:D-transform}
\E\bigl[e^{-\theta A}\bigr]
=\int_0^\infty e^{-t}\,e^{-\theta ct}\,e^{t\kappa(\theta)}\,dt
=\int_0^\infty e^{-t\,D(\theta)}\,dt=\frac{1}{D(\theta)},
\qquad
D(\theta):=1+\theta c-\kappa(\theta),
\end{equation}
whenever $\kappa(\theta)<\infty$ and $D(\theta)>0$; the integral diverges
otherwise. Since $\kappa$ is convex, \emph{$D$ is concave}, with $D(0)=1$ and
$D'(\theta)=c-\kappa'(\theta)$. Most of the analysis of this section is based on understanding properties of this function.

Indeed, both conditions of Proposition~\ref{prop:linear-existence} are
statements about $D$. Condition \eqref{eq:cond1}, namely
$\E[e^{-\gamma A}]=\tfrac12$, is
\begin{equation}\label{eq:Dgamma}
D(\gamma)=2 ,
\end{equation}
and differentiating \eqref{eq:D-transform} gives
$\E[Ae^{-\gamma A}]=D'(\gamma)/D(\gamma)^2$, so that the supercritical
condition \eqref{eq:cond2} is $D'(\gamma)>0$ and the critical condition
\eqref{eq:cond2ct} is $D'(\gamma)=0$. Combining this with
Lemma~\ref{lemma:deq_main} we obtain the following, which is all that
Sections~\ref{sec:brownian} and~\ref{sec:BS2024} will use.

\begin{proposition}[Travelling waves for a general L\'evy driver]
\label{prop:summary}
A travelling wave with decay parameter $\gamma>0$ and speed $c$, of the form
\begin{equation}\label{eq:sync-profile}
h(x)=\Prob(\chi>x)=1-\E\bigl[e^{-We^{-\gamma x}}\bigr],
\end{equation}
exists exactly when $D(\gamma)=2$. It is supercritical when $D'(\gamma)>0$,
in which case $W$ is the unique solution of \eqref{eq:deq_main_aux} with
$\E[W]<\infty$, and critical when $D'(\gamma)=0$, in which case $W$ is the
unique solution with $\E[W]=\infty$.
\end{proposition}

\medskip
Solving \eqref{eq:Dgamma} for $c$ yields the \emph{speed--decay relation}
\begin{equation}\label{eq:BS_speed}
c=v(\gamma):=\frac{1+\kappa(\gamma)}{\gamma},
\end{equation}
along which $v'(\gamma)=-D'(\gamma)/\gamma$. The supercritical regime is
therefore the decreasing branch of $v$, the critical regime its minimiser
$\gamma^*$, and $c^*:=v(\gamma^*)$ is the minimal speed; for $c<c^*$ the
relation \eqref{eq:BS_speed} has no solution and no wave of the form
\eqref{eq:sync-profile} exists.

\begin{figure}[ht]
\centering
\begin{tikzpicture}[xscale=2.2,yscale=1.25,>=stealth]
  \draw[->] (-1.85,0) -- (2.45,0) node[below] {$\theta$};
  \draw[->] (0,-0.9) -- (0,2.6);
  \draw[very thick,domain=-1.70:2.30,smooth,variable=\x]
       plot ({\x},{1+0.9*\x-0.13*\x*\x});
  \draw[dashed] (0,2) -- (1.39,2);
  \draw[dashed] (1.39,0) -- (1.39,2);
  \draw[dashed] (-1.39,0) -- (-1.39,-0.502);
  \draw[gray] (-1.39,-0.502) -- (0,1);
  \draw[gray] (0,1) -- (1.39,2);
  \node[circle,fill,inner sep=1pt] at (0,1) {};
  \node[circle,fill,inner sep=1pt] at (1.39,2) {};
  \node[circle,fill,inner sep=1pt] at (-1.39,-0.502) {};
  \node[circle,fill,inner sep=1pt] at (-0.974,0) {};
  \node[left] at (0,1) {$1$};
  \node[left] at (0,2) {$2$};
  \node[below] at (1.39,0) {$\gamma$};
  \node[below] at (-1.39,0) {$-\gamma$};
  \node[above left] at (-0.974,0) {$-\mu^*$};
  \node[right] at (2.30,{1+0.9*2.30-0.13*2.30*2.30}) {$D$};
\end{tikzpicture}
\caption{The concave function $D(\theta)=1+\theta c-\kappa(\theta)$, with
$D(0)=1$. A travelling wave of decay parameter $\gamma$ requires
$D(\gamma)=2$, and is supercritical or critical according as the curve is
still rising at $\gamma$ or has become flat there. The opposite tail is
governed by $-\mu^*$, the point at which $D$ first vanishes to the left of
the origin. The two chords have equal horizontal length $\gamma$, and
concavity forces the left one to fall at least as steeply as the right one
rises; hence $D(-\gamma)\le0$ and $\mu^*\le\gamma$
(Lemma~\ref{lem:mu-star}).}
\label{fig:D}
\end{figure}

\begin{remark}[Branching representation of $W$]\label{rem:BRW}
Proposition~\ref{prop:linear-existence} asserts the existence of $W$; the
following provides a way to construct it and in addition provides a connection to the
branching literature. Consider the branching L\'evy process in which each
particle lives an $\mathrm{Exp}(1)$ lifetime, moves during its life as
$\Gamma$, and at death produces two offspring at its position; let
$S_1(t),\dots,S_{\mathcal N_t}(t)$ be the positions of the particles alive at
time $t$. A standard calculation (either direct or using many-to-one technique) gives
$\E\bigl[\sum_ie^{\gamma S_i(t)}\bigr]=e^{(1+\kappa(\gamma))t}$, so that
\[
M_t=e^{-\gamma c\,t}\sum_{i=1}^{\mathcal N_t}e^{\gamma S_i(t)}
\]
is a martingale precisely when $\gamma c=1+\kappa(\gamma)$, i.e.\ when
$D(\gamma)=2$. Conditioning on the first branching event shows that its
almost sure limit satisfies \eqref{eq:deq_main_aux}. In the supercritical case
the tilted reproduction number $\widehat N=2e^{-\gamma A}$ has mean $1$ and,
provided $\E[A^-e^{-\gamma A}]<\infty$ --- which follows from
$\log^+(2e^{-\gamma A})\le\log2+\gamma A^-$ --- satisfies the
Kesten--Stigum criterion, so that $W>0$ almost surely on survival with
$\E[W]<\infty$ \cite{Biggins1}. In the critical case $M_t\to0$ and one uses
instead the derivative martingale
$\widetilde M_t=e^{-\gamma c\,t}\sum_i(ct-S_i(t))e^{\gamma S_i(t)}$, whose
limit is positive on survival under the extra moment condition of
\eqref{eq:cond2ct} \cite{Biggins2}, with infinite mean.
\end{remark}

\subsection{Tails of the travelling waves}\label{sec:mustar}

The representation \eqref{eq:sync-profile} reads the two tails of the wave off
the two ends of the Laplace transform $\varphi_W(s)=\E[e^{-sW}]$: the
behaviour at $s\to0$ governs $h(x)$ as $x\to+\infty$, and gives
$h(x)\sim\E[W]e^{-\gamma x}$ whenever $\E[W]<\infty$; the behaviour at
$s\to\infty$ governs $1-h(x)$ as $x\to-\infty$. The second is \emph{not}
governed by $\gamma$, and the relevant exponent is again visible in
Figure~\ref{fig:D}.

\begin{lemma}\label{lem:mu-star}
Let $\gamma>0$ satisfy $D(\gamma)=2$, let $h$ be as in \eqref{eq:sync-profile},
and set $\mu^*:=\sup\{\mu\ge0:\E[e^{\mu A}]<\infty\}$. Write $\eta:=-\chi$, so
that \eqref{eq:main_deq_kpp} reads $\eta\overset{d}{=}U+A$ with
$U:=\min\{\eta_1,\eta_2\}$ independent of $A$.
\begin{enumerate}
\item[(a)] If $\E[W]<\infty$ then $h(x)\sim\E[W]\,e^{-\gamma x}$ as
$x\to+\infty$. More generally $h(x)\le\E[W^{\alpha}]\,e^{-\alpha\gamma x}$ for
every $\alpha\in(0,1]$, so that $h$ decays at the exponential rate
$\alpha\gamma$ as soon as $\E[W^{\alpha}]<\infty$ --- which holds for every
$\alpha<1$ whenever $\E[W]<\infty$.
\item[(b)] Assume $\mu^*>0$, and, for the last assertion, that
$\E[W^{\alpha}]<\infty$ for some $\alpha\in(0,1)$, which by \emph{(a)} is
automatic when $\E[W]<\infty$. Then $\E[e^{\beta\eta}]<\infty$ for $\beta<\mu^*$
and $\E[e^{\beta\eta}]=\infty$ for $\beta>\mu^*$, while $\E[e^{\beta U}]<\infty$
for every $\beta<2\mu^*$. If moreover $\Prob(A>t)\sim L\,t^{k}e^{-\mu^*t}$ as
$t\to\infty$, for some $L>0$ and $k\ge0$, then
\[
1-h(x)\;\sim\;L\,\E\bigl[e^{\mu^*U}\bigr]\,|x|^{k}e^{-\mu^*|x|}
\qquad(x\to-\infty).
\]
\item[(c)] $\mu^*\le\gamma$, with equality if and only if $\kappa$ is affine.
\end{enumerate}
\end{lemma}

\begin{proof}
(a) Define $s: x \to e^{-\gamma x}$. Dominated convergence applied to
$h(x)/s(x)=\E[W(1-e^{-Ws(x)})/(Ws(x))]$, the ratio lying in $[0,1]$ and tending to $1$ as $x \to +\infty$,
gives $\E[W]$. For the second assertion, we use $1-e^{-u}\le u^{\alpha}$ for
$u\ge0$ and $\alpha\in(0,1]$, whence $h(x)\le\E[W^{\alpha}]e^{-\alpha\gamma x}$;
and $\E[W^{\alpha}]<\infty$ for every $\alpha<1$ when $\E[W]<\infty$, by
Jensen's inequality.

(b) Put $m(\beta):=\E[e^{\beta\eta}]$ and
$\beta^*:=\sup\{\beta\ge0:m(\beta)<\infty\}$. Independence gives
$m(\beta)=\E[e^{\beta U}]\E[e^{\beta A}]$ with both factors positive, so
$\beta^*\le\mu^*$; and $\min(y_1,y_2)\le\sqrt{y_1y_2}$ with
$y_i=e^{\beta\eta_i}$ gives $\E[e^{\beta U}]\le m(\beta/2)^2$, whence
\begin{equation}\label{eq:tail-doubling}
m(\beta)\le m(\beta/2)^2\,\E[e^{\beta A}].
\end{equation}
So $m$ is finite on $[0,\min\{2\beta^*,\mu^*\})$, i.e.
$\min\{2\beta^*,\mu^*\}\le\beta^*$; if $\beta^*>0$ this forces
$\beta^*=\mu^*$, and \eqref{eq:tail-doubling} then bounds $\E[e^{\beta U}]$ for
$\beta<2\mu^*$.

To prove $\beta^*>0$, fix $a\in(0,\mu^*)$ and $C\ge1$ with
$\Prob(A>t)\le Ce^{-a t}$, and write $\bar h(x):=\Prob(\eta>x)$. Splitting
$\{U+A>x\}$ at $x/2$ gives $\bar h(x)\le\bar h(x/2)^2+Ce^{-a x/2}$. Choose
$x_0$ with $\bar h(x_0)\le\tfrac14$ and $Ce^{-a x_0/2}\le\tfrac14$, and set
$\beta:=\min\{a/4,(\log2)/x_0\}$; then $\bar h(x_0)\le\tfrac12e^{-\beta x_0}$
and the bound propagates along $x_n=2^nx_0$, since
$\tfrac14e^{-2\beta x_n}+Ce^{-a x_n}\le\tfrac12e^{-\beta x_{n+1}}$. By
monotonicity $\bar h(x)\le\tfrac12e^{-\beta x/2}$ for $x\ge x_0$, so
$\beta^*\ge\beta/2>0$.

For the asymptotic, $\phi(t):=e^{\mu^*t}t^{-k}\Prob(A>t)\to L$ and, for $x\ge1$,
\[
e^{\mu^*x}x^{-k}\,\bar h(x)
=\E\bigl[e^{\mu^*U}\phi(x-U)\bigl((x-U)/x\bigr)^{k}\bigr],
\]
whose integrand tends to $L\,e^{\mu^*U}$ and is dominated by
$C\,e^{\mu^*U}(2+|U|)^{k}$, using
$\Prob(A>t)\le C(1+t^+)^{k}e^{-\mu^*t}$. This is integrable: on $\{U\ge0\}$
because $\E[e^{\beta U}]<\infty$ for some $\beta\in(\mu^*,2\mu^*)$, and on
$\{U<0\}$ because $\Prob(U<-t)\le2h(t)$ decays exponentially by (a). Since
$1-h(x)=\bar h(-x)$, this is the claim.

(c) As $D(0)=1$ and $D(\gamma)=2$, concavity implies
$\bigl(1-D(-\gamma)\bigr)/\gamma\ge1/\gamma$; hence $D(-\gamma)\le0$,
$\E[e^{\gamma A}]=\infty$ and $\mu^*\le\gamma$. Equality forces $D(-\gamma)=0$
and $D$ affine on $[-\gamma,\gamma]$, i.e.\ $\kappa$ affine.
\end{proof}

\begin{remark}\label{rem:left-tail}
The exponent $\mu^*$ is the one produced by linearising the travelling-wave
equation at the stable state $h=1$. Indeed $F:=1-h$ satisfies
$F(x)=\E[F(x+A)^2]$ and, $T$ being exponential, the resolvent identity
$\E[G(x+cT-\Gamma(T))]=(I-\mathcal A)^{-1}G(x)$, with
$\mathcal Af=cf'+\tfrac{\sigma^2}{2}f''$ for Brownian $\Gamma$, turns this
into
\begin{equation}\label{eq:tw-ode}
F-cF'-\tfrac{\sigma^2}{2}F''=F^2,
\qquad\text{equivalently}\qquad
\tfrac{\sigma^2}{2}h''+ch'+h-h^2=0,
\end{equation}
the travelling-wave equation of \eqref{KPP} --- a useful consistency check on
the approach. As $F\to0$ at $-\infty$ the quadratic term is negligible there
and $F$ decays like $e^{\mu x}$ with $\tfrac{\sigma^2}{2}\mu^2+c\mu-1=0$,
whose positive root $\mu^*=(\sqrt{c^2+2\sigma^2}-c)/\sigma^2$ is exactly the
exponent of Lemma~\ref{lem:mu-star}: here $\kappa(-\mu)=\sigma^2\mu^2/2$ is
finite for every $\mu\ge0$, so that by \eqref{eq:D-transform} the expectation
$\E[e^{\mu A}]$ is finite precisely while
$D(-\mu)=1-\mu c-\sigma^2\mu^2/2>0$. Classical F-KPP results
\cite{AW,U} then give $1-h(x)\sim Ce^{\mu^*x}$. For a general L\'evy driver
\eqref{eq:tw-ode} becomes a renewal equation and we expect the same
conclusion by the exponential change of measure of the renewal theorem.
\end{remark}

\subsection{Brownian motion}\label{sec:brownian}

For a Brownian motion of variance $\sigma^2$ one has
$\kappa(\theta)=\sigma^2\theta^2/2$, so that \eqref{eq:BS_speed} reads
\[
c=v(\gamma)=\frac1\gamma+\frac{\sigma^2\gamma}{2}.
\]
This is strictly convex on $(0,\infty)$ and minimised at
$\gamma^*=\sqrt2/\sigma$, with $c^*=v(\gamma^*)=\sqrt2\,\sigma$: the classical
minimal speed of \eqref{KPP}, which also confirms the normalisation of
\eqref{eq:generator}. For $c>c^*$ the wave decays like $e^{-\gamma x}$ with
$\gamma\in(0,\gamma^*)$ the smaller root of $v(\gamma)=c$, and at $c=c^*$ like
$xe^{-\gamma^*x}$. The opposite tail decays at rate
$\mu^*=(\sqrt{c^2+2\sigma^2}-c)/\sigma^2<\gamma$, by
Lemma~\ref{lem:mu-star} and Remark~\ref{rem:left-tail}.

\subsection{Compound Poisson: the Baryshnikov--Stolyar model}\label{sec:BS2024}

Consider the model of Example~\ref{ex:mixed}, in which each particle jumps at
rate $\lambda>0$ by an independent amount with law $F_X$, $\E[X]<\infty$, and
at rate $1$ copies the position of a uniformly chosen other particle whenever
that position exceeds its own. The latter is a probabilistic variant of the model of \cite{GJM}, which converges to the F-KPP equation.
Define:
\[
\kappa(\theta)=\lambda\bigl(\E[e^{\theta X}]-1\bigr),
\]
finite for $\theta<\gamma_0:=\sup\{\gamma:\E[e^{\gamma X}]<\infty\}$, and
\eqref{eq:BS_speed} becomes
\begin{equation}\label{eq:BS_gamma_eq}
c=v(\gamma)=\frac{1+\lambda\bigl(\E[e^{\gamma X}]-1\bigr)}{\gamma}.
\end{equation}

Everything below depends on the minimum of $v$ being attained \emph{inside}
$(0,\gamma_0)$ rather than at its right endpoint. We impose this as a
hypothesis on the jump law, stated in the explicit form in which it is
checked.

\medskip
\noindent\textbf{Assumption I.}
\emph{$X\ge0$ with $\Prob(X>0)>0$, $\gamma_0>0$, and}
\begin{equation}\label{eq:A}
\lambda\,\E\bigl[Xe^{\gamma X}\bigr]\;>\;v(\gamma)
\qquad\text{for some }\gamma\in(0,\gamma_0).
\end{equation}

\medskip
The following lemma says what \eqref{eq:A} is worth: it is precisely the
condition under which the minimal speed is produced by an interior decay
parameter, and it is implied by the properties of $X$ that one checks in
practice.

\begin{lemma}\label{lem:steep}
Let $X\ge0$ with $\Prob(X>0)>0$ and $\gamma_0>0$. The function $v$ of
\eqref{eq:BS_gamma_eq} is strictly convex on $(0,\gamma_0)$, with
$v(0^+)=+\infty$ and $v'(\gamma)=-D'(\gamma)/\gamma$. Moreover the following
are equivalent:
\begin{enumerate}
\item[\emph{(i)}] \eqref{eq:A} holds;
\item[\emph{(ii)}] $v$ attains its infimum over $(0,\gamma_0)$ at a
(necessarily unique) interior point $\gamma^*$, characterised by
$v(\gamma^*)=\lambda\E[Xe^{\gamma^*X}]$, equivalently by $D'(\gamma^*)=0$;
\item[\emph{(iii)}] $v$ is not non-increasing on $(0,\gamma_0)$.
\end{enumerate}
Each of the following is sufficient for \eqref{eq:A}: $\gamma_0=\infty$,
which holds in particular for bounded and for Gaussian-tailed jumps;
$\E[e^{\gamma X}]\to\infty$ as $\gamma\uparrow\gamma_0$; or
$\E[Xe^{\gamma X}]\to\infty$ as $\gamma\uparrow\gamma_0$.
\end{lemma}

\begin{proof}
Write $\Lambda(\gamma):=\E[e^{\gamma X}]$ and
$g(\gamma):=\gamma^2\Lambda''(\gamma)-2\gamma\Lambda'(\gamma)+2\Lambda(\gamma)-2$,
so that $g(0)=0$ and
$g'(\gamma)=\gamma^2\Lambda'''(\gamma)=\gamma^2\E[X^3e^{\gamma X}]\ge0$; hence
$g\ge0$ and $v''(\gamma)=2\gamma^{-3}+\lambda\gamma^{-3}g(\gamma)>0$. Thus
$v'$ is strictly increasing, and $v'(\gamma)\sim-\gamma^{-2}\to-\infty$ as
$\gamma\downarrow0$. Since
$\gamma^2v'(\gamma)=\lambda\gamma\E[Xe^{\gamma X}]-\gamma\,v(\gamma)$,
condition \eqref{eq:A} says exactly that $v'>0$ somewhere, which by strict
monotonicity of $v'$ is equivalent to the existence of a unique zero of $v'$
in $(0,\gamma_0)$, and to (iii). The sufficient conditions follow since
$\lambda\gamma\E[Xe^{\gamma X}]-\gamma v(\gamma)\to\infty$ in each case.
\end{proof}

Under Assumption I we write $\gamma^*$ for the interior minimiser supplied by
Lemma~\ref{lem:steep}(ii), so that
\begin{equation}\label{eq:BS_gamma_star}
v(\gamma^*)=\lambda\,\E\bigl[Xe^{\gamma^*X}\bigr],
\end{equation}
and $c^*:=v(\gamma^*)$. When \eqref{eq:A} fails, $v$ decreases on all of
$(0,\gamma_0)$, the infimum $c^*=v(\gamma_0^-)$ is approached only at the
edge, \eqref{eq:BS_gamma_star} has no solution, and the front is \emph{pushed}
rather than pulled; that case is outside our analysis. Written out, the
supercritical condition $D'(\gamma)>0$ reads
\begin{equation}\label{eq:BS_supercrit}
\frac{1+\lambda\bigl(\E[e^{\gamma X}]-1\bigr)}{\gamma}
>\lambda\,\E\bigl[Xe^{\gamma X}\bigr],
\end{equation}
which holds precisely when $\gamma<\gamma^*$.

The statement below collects, for this model, the consequences of
Proposition~\ref{prop:summary}, Remark~\ref{rem:BRW} and
Lemma~\ref{lem:mu-star}.

\begin{proposition}[Existence, uniqueness, and representation of travelling waves]
\label{thm:BS_TW}
Let Assumption I hold, and write
\begin{equation}\label{eq:BS_TW_form}
h(x)=\Prob(\chi>x)=1-\E\bigl[e^{-Ve^{-\gamma x}}\bigr],
\end{equation}
$V$ being a solution of
\begin{equation}\label{eq:BS_linear}
V \overset{d}{=} (V_1+V_2)\,e^{-\gamma A},
\qquad A=cT-\Gamma(T).
\end{equation}
\begin{enumerate}
\item[\emph{(i)}] \textbf{Supercritical regime} $(c>c^*)$. There is a unique
$\gamma\in(0,\gamma^*)$ with $v(\gamma)=c$, and the additive martingale of
Remark~\ref{rem:BRW} converges almost surely to a limit $W>0$ with
$\E[W]<\infty$. The profile is
\begin{equation}\label{eq:BS_TW_super}
h(x)=1-\E\bigl[e^{-We^{-\gamma x}}\bigr],
\end{equation}
and $W$ is the unique solution of \eqref{eq:BS_linear} with finite mean, up
to a multiplicative constant.

\item[\emph{(ii)}] \textbf{Critical regime} $(c=c^*)$. The additive
martingale converges to zero. If $\E[|A|^{1+\delta}e^{-\gamma^*A}]<\infty$ for
some $\delta>0$, the derivative martingale of Remark~\ref{rem:BRW} converges
almost surely to $\widetilde W>0$ with $\E[\widetilde W]=\infty$, and
\begin{equation}\label{eq:BS_TW_crit}
h(x)=1-\E\bigl[e^{-\widetilde We^{-\gamma^*x}}\bigr].
\end{equation}

\item[\emph{(iii)}] \textbf{No wave below $c^*$.} For $c<c^*$ the relation
\eqref{eq:BS_speed} has no solution $\gamma>0$ and no wave of the form
\eqref{eq:BS_TW_form} exists.

\item[\emph{(iv)}] \textbf{Tails.} In the supercritical regime
$h(x)\sim\E[W]e^{-\gamma x}$ as $x\to+\infty$, and in the critical regime
$h(x)\sim Cxe^{-\gamma^*x}$ for some $C>0$. The opposite tail decays at the
rate $\mu^*$ of Lemma~\ref{lem:mu-star}, and $\mu^*<\gamma$ strictly, since
$\kappa$ is affine only if $X\equiv0$; see also Remark~\ref{rem:left-tail}.
\end{enumerate}
\end{proposition}

\begin{proof}
Taking expectations in \eqref{eq:BS_linear}, assuming $\E[V]\in(0,\infty)$,
gives $\E[V]=2\E[V]\E[e^{-\gamma A}]$ and hence $\E[e^{-\gamma A}]=\tfrac12$,
which by \eqref{eq:D-transform} is $D(\gamma)=2$, i.e.\ \eqref{eq:BS_speed}.
Parts (i) and (ii) are then Proposition~\ref{prop:summary} together with
Remark~\ref{rem:BRW}, the supercritical condition $D'(\gamma)>0$ being
\eqref{eq:BS_supercrit}, valid exactly for $\gamma<\gamma^*$, and the critical
condition $D'(\gamma^*)=0$ being \eqref{eq:BS_gamma_star}. For $c>c^*$ the
equation $v(\gamma)=c$ has, by Lemma~\ref{lem:steep}, exactly one root
$\gamma_-\in(0,\gamma^*)$ and at most one root $\gamma_+\in(\gamma^*,\gamma_0)$.
At $\gamma_+$ one has $D'(\gamma_+)<0$, so that $\E[Ae^{-\gamma_+A}]<0$ and
neither condition of Proposition~\ref{prop:linear-existence} holds:
\eqref{eq:BS_linear} has no solution of finite mean there. It does have one of
infinite mean, but that solution yields no new wave. Indeed
$\theta:=\gamma_-/\gamma_+\in(0,1)$ satisfies
$2\,\E[e^{-\theta\gamma_+A}]=2/D(\gamma_-)=1$ while
$\tfrac{d}{d\theta}2\E[e^{-\theta\gamma_+A}]<0$ there, since $D'(\gamma_-)>0$;
the associated fixed point therefore has a regularly varying tail of index
$\theta$ \cite{Liu2001}, so that
$1-\E\bigl[e^{-Ve^{-\gamma_+x}}\bigr]$ is of order
$e^{-\theta\gamma_+x}=e^{-\gamma_-x}$ as $x\to+\infty$: the profile obtained at
$\gamma_+$ is the one already obtained at $\gamma_-$. Part (iii) is the strict
convexity of $v$ (Lemma~\ref{lem:steep}). In part (iv), the right tail follows from
\eqref{eq:BS_TW_super} by writing $s=e^{-\gamma x}\to0^+$ and applying
dominated convergence to
$h(x)/e^{-\gamma x}=\E[W(1-e^{-Ws})/(Ws)]\to\E[W]$; the critical tail is a
standard consequence of the derivative-martingale construction and the
infinite mean of $\widetilde W$, see \cite{Biggins2}; and the statements about
the left tail are Lemma~\ref{lem:mu-star} applied to $A=cT-\Gamma(T)$, whose
$\kappa(\theta)=\lambda(\E[e^{\theta X}]-1)$ is affine only if $X\equiv0$.
\end{proof}

\begin{remark}[Uniqueness of the profile]\label{rem:BS_uniqueness}
For $c>c^*$, uniqueness of the wave up to spatial shifts follows from the
uniqueness of $V$ in Proposition~\ref{prop:linear-existence}: multiplicative
constants for $V$ correspond precisely to translations of $h$. At $c=c^*$ the
analogous statement holds through the derivative martingale; see
\cite[Proposition~2.2(ii)]{liu1999asymptotic}.
\end{remark}

\begin{remark}\label{rem:BS_relation}
The model interpolates between pure F-KPP dynamics ($\lambda=0$) and a regime
dominated by independent jumps ($\lambda\to\infty$); the first case is treated
in Example~\ref{ex:pure_copying}. In general the wave retains the exponential
structure \eqref{eq:BS_TW_super}, with both the speed and the shape parameter
modulated by the jump component through \eqref{eq:BS_speed}; the branching
construction of Remark~\ref{rem:BRW} also gives a practical way of simulating
the profile through the law of $W$.
\end{remark}


\begin{example}[Pure copying, $\lambda=0$]\label{ex:pure_copying}
Here $\Gamma\equiv0$, so $\kappa\equiv0$, $D(\theta)=1+\theta c$ and
\eqref{eq:Dgamma} reads $\gamma c=1$. Every $\gamma>0$ is admissible, with
speed $c=1/\gamma$: a one-parameter family of travelling waves, which is the
classical F-KPP family without diffusion. Selection of a member depends on the
initial condition \cite{AW}. Note that $D$ is affine here, the equality case
of Lemma~\ref{lem:mu-star}(c): indeed $\mu^*=1/c=\gamma$, and this is the only
situation in which the two tails share an exponent.
\end{example}


\begin{example}[Exponential jumps]
For $X\sim\mathrm{Exp}(1)$ one has $\E[e^{\gamma X}]=1/(1-\gamma)$ for
$\gamma<1$, so that \eqref{eq:BS_speed} gives
\[
v(\gamma)=\frac1\gamma+\frac{\lambda}{1-\gamma},
\]
whose minimiser and minimum are explicit. This recovers the results of
\cite{baryshnikov2025large}; for $\lambda\ll1$, expanding
$v(\gamma)=\gamma^{-1}+\lambda\gamma^{-1}(\E[e^{\gamma X}]-1)+O(\lambda^2)$
gives the perturbation of the F-KPP speed caused by a small jump activity, for
a general jump law.
\end{example}

\section{Travelling Waves for the Power-of-2 Model} \label{sec:power-of-2}

\subsection{Characterisation through Distributional Equations}

Consider the power-of-2 model described in Example~\ref{ex:power2}, the
particles being driven by a L\'evy process $\Gamma$ with
$\E|\Gamma(1)|<\infty$. Brownian motion is the main example, and the one for
which the mean-field limit has been established.
By Proposition~\ref{prop:TW-FP} (with $\alpha=1$), a travelling wave profile
(i.e., the stationary distribution of a tagged particle in the frame moving
with speed $c$) is exactly a solution of the distributional equation:

\begin{equation} \label{eq:distr1}
\chi \overset{d}{=} \chi_1 + \I{\chi_1 \le \chi_2} X - cT + \Gamma(T),
\end{equation}
where:
\begin{itemize}
\item $\chi_1, \chi_2, X, T$ are independent random variables,
\item $\Gamma(T)$ is the displacement of the driving process over a time
interval of length $T$ (a Brownian motion with variance $\sigma^2$ per unit
time in the main example),
\item $X$ has distribution $F$ with $\E[X]=1$,
\item $\chi_1, \chi_2, \chi$ have the same distribution $h$,
\item $T \sim \mathrm{Exp}(1)$ is the time between the interactions initiated by the tagged particle.
\end{itemize}

A crucial simplification occurs by noting that for this model the speed is
determined by the mechanism alone. Under the normalisation of
\eqref{eq:generator} it equals
\begin{equation}\label{eq:c-pow2}
c=\gamma_1+\bar\beta=\gamma_1+\frac{\E[X]}{2},\qquad\gamma_1:=\E\,\Gamma(1),
\end{equation}
by \eqref{eq:speed-pow2} and \eqref{eq:pow2-speed}, which is $\E[X]/2=1/2$ for
a centred driver such as driftless Brownian motion; this is the analogue, in
our normalisation, of the value $c=\E[X]$ obtained in \cite{Stolyar2023} from
conservation of mass, which corresponds to interaction clocks of rate $2$
(see Remark~\ref{rem:normalisation}). As shown in Step~3 of the proof of
Proposition~\ref{prop:equiv_power2}, the value \eqref{eq:c-pow2} is in fact
forced by the balance equation of the tagged particle itself, independently of
Lemma~\ref{lem:bary-gen}. With interaction clocks of rate $\alpha$ that same
computation gives $c=\gamma_1+\alpha\E[X]/2$; since the functionals $M$ and
$\overline H$ appearing below depend on the normalisation only through $c$,
the decomposition of Proposition~\ref{prop:equiv_power2} keeps its form.

The statement below involves two functionals of the driving process seen in
the wave frame, $Y(t):=\Gamma(t)-ct$: its all-time supremum
$M=\sup_{t\ge0}Y(t)$, and the equilibrium law $\overline H\ge0$ of the
descending ladder height of $Y$, defined in \eqref{eq:Hbar} of
Appendix~\ref{app:ladder}. All we use about them is
Lemma~\ref{lem:adjoint}, together with the fact that $\overline H\equiv0$ if
and only if $Y$ has no negative jumps: a process reaching its new infima
continuously has a ladder height which is a pure drift.

We say that two distributional equations are \emph{equivalent} if there exists
a bijective transformation mapping solutions of one equation to solutions of
the other. The following proposition provides a rather surprising way to
connect travelling waves in the models with synchronisation from the previous
section and the previously described power-of-2 interaction.

\begin{proposition}[Equivalent formulation, general L\'evy driver] \label{prop:equiv_power2}
Let $\Gamma$ be a L\'evy process with $\E|\Gamma(1)|<\infty$, write
$\gamma_1:=\E\,\Gamma(1)$ and $Y(t):=\Gamma(t)-ct$, and assume that the law of
$\chi$ is atomless (as is the case as soon as $\Gamma$ has a non-degenerate
Gaussian component); $h(x)=\Prob(\chi>x)$ denotes its survival function. Then
the distributional equation \eqref{eq:distr1} admits a solution only if
\begin{equation}\label{eq:c-pow2-gen}
c=\gamma_1+\bar\beta,\qquad \bar\beta=\frac{\E[X]}{2}=:\mu ,
\end{equation}
and in that case it is equivalent to
\begin{equation} \label{eq:distr2-gen}
\chi-\overline H\;\overset{d}{=}\;\min(\chi_1, \chi_2) + \overline{X} + M,
\end{equation}
where $\min(\chi_1,\chi_2)$, $\overline X$, $\overline H$ and $M$ are mutually
independent and:
\begin{itemize}
\item $\overline{X}$ is the integrated tail (equilibrium) random variable with survival function
$$
\Prob(\overline{X} > x) = \frac{1}{\E[X]}\int_x^\infty \overline{F}(y)dy, \quad \text{where } \overline{F}(x) = 1 - F(x);
$$
\item $M$ is the all-time supremum of the driving process in the wave frame,
\begin{equation}\label{eq:M-def}
M:=\sup_{t\ge0}\bigl(\Gamma(t)-ct\bigr) ,
\end{equation}
almost surely finite because $\E\,Y(1)=-\mu<0$ by \eqref{eq:c-pow2-gen};
\item $\overline H$ is the equilibrium descending ladder height of $Y$, as in
\eqref{eq:Hbar}.
\end{itemize}
\end{proposition}

\begin{corollary}[Spectrally positive drivers]\label{cor:equiv_pow2_sp}
If in addition $\Gamma$ has no negative jumps, then $\overline H\equiv0$ and
\eqref{eq:distr2-gen} reduces to
\begin{equation} \label{eq:distr2}
\chi \overset{d}{=} \min(\chi_1, \chi_2) + \overline{X} + M ,
\end{equation}
the law of $M$ being given by the generalised Pollaczek--Khinchine formula
\begin{equation}\label{eq:PK}
\E\bigl[e^{-sM}\bigr]=\frac{\mu\,s}{sc+\kappa(-s)},\qquad s>0,
\qquad \kappa(\theta):=\log\E\bigl[e^{\theta\Gamma(1)}\bigr],
\end{equation}
$\kappa$ being finite on $\theta\le0$. In particular, when $\Gamma$ is a
driftless Brownian motion of variance $\sigma^2$ one has $\gamma_1=0$,
$c=\E[X]/2$, $\kappa(-s)=\sigma^2s^2/2$, and \eqref{eq:PK} reduces to
$\E[e^{-sM}]=\bigl(1+\sigma^2s/(2c)\bigr)^{-1}$: the variable $M$ is then
exponentially distributed with mean $\sigma^2/(2c)=\sigma^2/\E[X]$, which is
$\sigma^2$ under the normalisation $\E[X]=1$ of Example~\ref{ex:power2}.
\end{corollary}

\begin{proof}
Only \eqref{eq:PK} requires an argument, and it is the generalised
Pollaczek--Khinchine formula for the stationary content of a spectrally
positive L\'evy process reflected at its infimum (for instance, see \cite[Ch.~IX, Corr. 3.4]{Asmussen03}).
It is
also what Lemma~\ref{lem:adjoint} gives formally on the exponential test
functions $f=e^{-s\,\cdot}$, $s>0$, which do not belong to $\mathcal T$: with $\overline H\equiv0$ and
$\mathcal A_Ye^{-sx}=\psi_Y(s)e^{-sx}$, where
$\psi_Y(s)=\log\E[e^{-sY(1)}]=sc+\kappa(-s)$ is finite for $s\ge0$ because $Y$
has no negative jumps, \eqref{eq:adjoint} reads
$\psi_Y(s)\E[e^{-sM}]=\mu s$.
\end{proof}

\begin{proof}[Proof of Proposition~\ref{prop:equiv_power2}]
Write $U:=\min\{\chi_1,\chi_2\}$ and let $\mathcal A_Y$ be the generator of $Y$.
Throughout, $f$ ranges over the class $\mathcal T$ of \eqref{eq:testclass};
the two exceptions, $f=\mathrm{id}$ in Step~3 and $f=e^{i\lambda\,\cdot}$ in
Remark~\ref{rem:transform-form}, are obtained by localisation and by
approximation respectively.

\smallskip
\emph{Step 1.} For every
bounded measurable $g$, since $h(\chi_1)=\Prob(\chi_2>\chi_1\mid\chi_1)$ and
$\chi_1=U$ on $\{\chi_1<\chi_2\}$,
\begin{equation}\label{eq:min-density}
\E\bigl[g(\chi)h(\chi)\bigr]
=\E\bigl[g(\chi_1)\I{\chi_1<\chi_2}\bigr]
=\E\bigl[g(U)\I{\chi_1<\chi_2}\bigr]
=\tfrac12\,\E\bigl[g(U)\bigr],
\end{equation}
the last equality because exchanging $\chi_1$ and $\chi_2$ shows that
$\{\chi_1<\chi_2\}$ and $\{\chi_2<\chi_1\}$ contribute equally, and together they
exhaust the space since $\Prob(\chi_1=\chi_2)=0$. Equivalently,
$h(y)\,\Prob(\chi\in dy)=\tfrac12\,\Prob(U\in dy)$.

\smallskip
\emph{Step 2.} In the wave frame the tagged particle
moves as $Y$ and, at rate $1$, samples an independent $\chi_2\sim h$ and is
pushed up by $X$ whenever it lies below it; from the position $y$ it therefore
jumps at rate $h(y)$, and its generator is
\begin{equation}\label{eq:tagged-generator}
\mathcal Lf(y)=\mathcal A_Yf(y)+h(y)\,\E\bigl[f(y+X)-f(y)\bigr].
\end{equation}
This is $\mathcal A_Yf+K_\nu f-f$ in the notation of
Proposition~\ref{prop:TW-FP}, $\nu$ being the law of $\chi$; since every
$f\in\mathcal T$ lies in $C^2_0$ by the Riemann--Lebesgue lemma, the
equivalence (b)$\Leftrightarrow$(c) there shows that \eqref{eq:distr1} says
precisely that $h$ is invariant for \eqref{eq:tagged-generator}, i.e.\ that
$\E_h[\mathcal Lf(\chi)]=0$. By Step~1 this reads
\begin{equation}\label{eq:pow2-mid}
\E\bigl[\mathcal A_Yf(\chi)\bigr]+\tfrac12\,\E\bigl[f(U+X)-f(U)\bigr]=0 .
\end{equation}

\smallskip
\emph{Step 3.} Let $f$ be absolutely
continuous with bounded derivative. For $X\ge0$ with $\E[X]<\infty$, writing
$f(y+X)-f(y)=\int_0^Xf'(y+u)\,du$ and applying Fubini gives, for every
$y\in\mathbb R$,
\begin{equation}\label{eq:equilibrium-jump}
\E\bigl[f(y+X)\bigr]-f(y)
=\int_0^\infty f'(y+u)\,\Prob(X>u)\,du
=\E[X]\,\E\bigl[f'(y+\overline X)\bigr],
\end{equation}
the last equality because $\overline X$ has density $\Prob(X>u)/\E[X]$ on
$(0,\infty)$. A jump therefore contributes to the balance an average of $f'$
over the levels it sweeps, and $\overline X$ is the law of a uniformly chosen
such level --- which is where the equilibrium distribution comes from.

Apply \eqref{eq:equilibrium-jump} at $y=U$ and average over $U$, which is
legitimate since $X$ is independent of $U$. Writing $S:=U+\overline X$, with
$\overline X$ independent of $U$, and recalling $\mu=\E[X]/2$ from
\eqref{eq:c-pow2-gen}, the interaction term of \eqref{eq:pow2-mid} becomes
\[
\tfrac12\,\E\bigl[f(U+X)-f(U)\bigr]
=\tfrac12\,\E[X]\;\E\bigl[f'(U+\overline X)\bigr]
=\mu\,\E\bigl[f'(S)\bigr] ,
\]
the factor $\tfrac12$ of the balance equation and the factor $\E[X]$ of
\eqref{eq:equilibrium-jump} combining into $\mu$. Hence \eqref{eq:pow2-mid}
reads
\begin{equation}\label{eq:pow2-star}
\E\bigl[\mathcal A_Yf(\chi)\bigr]=-\mu\,\E\bigl[f'(S)\bigr].
\end{equation}
Taking $f=\mathrm{id}$ in \eqref{eq:pow2-mid} (after the usual localisation)
gives $\E\,Y(1)+\E[X]/2=0$, that is $\gamma_1-c+\E[X]/2=0$, which is
\eqref{eq:c-pow2-gen}; in particular $\E\,Y(1)=-\mu<0$, so that $M<\infty$
almost surely and Lemma~\ref{lem:adjoint} applies to $Y$.

\smallskip
\emph{Step 4.} Let $\overline H$ be as in \eqref{eq:Hbar},
independent of everything else, and apply \eqref{eq:pow2-star} to
$\tilde f:=\E[f(\,\cdot-\overline H)]$. Since $\mathcal A_Y$ commutes with
translations, $\mathcal A_Y\tilde f(x)=\E[\mathcal A_Yf(x-\overline H)]$, so
that
\begin{equation}\label{eq:step4a}
\E\bigl[\mathcal A_Yf(\chi-\overline H)\bigr]
=-\mu\,\E\bigl[\tilde f'(S)\bigr]
=-\mu\,\E\bigl[f'(S-\overline H)\bigr].
\end{equation}
On the other hand, applying Lemma~\ref{lem:adjoint} to the translate
$f(s+\,\cdot)$ and averaging over an independent $s\sim S$ gives
\begin{equation}\label{eq:step4b}
\E\bigl[\mathcal A_Yf(S+M)\bigr]=-\mu\,\E\bigl[f'(S-\overline H)\bigr].
\end{equation}
Comparing \eqref{eq:step4a} and \eqref{eq:step4b},
$\E[\mathcal A_Yf(\chi-\overline H)]=\E[\mathcal A_Yf(S+M)]$ for every test
function $f$. Let $\rho$ be the difference of the two laws: it is a finite
signed measure with $\rho(\mathbb R)=0$ and $\int\mathcal A_Yf\,d\rho=0$ for
all $f$, i.e.\ $\hat\rho(\lambda)\Psi_Y(\lambda)=0$ for every $\lambda$. The
zero set of $\Psi_Y$ is contained in
$\{\lambda:\varphi_{Y(1)}(\lambda)=1\}$, a closed subgroup of $\mathbb R$
which is not all of $\mathbb R$, since $\varphi_{Y(1)}\equiv1$ would force
$Y(1)=0$ almost surely against $\E\,Y(1)=-\mu<0$; being a proper closed
subgroup it is discrete, so $\hat\rho$ vanishes off a discrete set and, being
continuous, vanishes identically. Thus $\rho=0$. This is
\eqref{eq:distr2-gen}. Every step above is an equivalence, so the two
formulations are equivalent.
\end{proof}

\begin{remark}[Transform form]\label{rem:transform-form}
Taking $f=e^{i\lambda\,\cdot}$ in the balance equation \eqref{eq:pow2-mid},
and using $\varphi_{\overline X}(\lambda)=(\varphi_X(\lambda)-1)/(i\lambda\E[X])$
together with \eqref{eq:min-density}, returns the transform identity
\begin{equation}\label{eq:pow2-final}
\varphi_\chi(\lambda)
=\underbrace{\frac{i\lambda\,\mu}{i\lambda c-\Psi(\lambda)}}_{=:R(\lambda)}
\cdot\varphi_{\overline X}(\lambda)\cdot\varphi_{\min\{\chi_1,\chi_2\}}(\lambda),
\end{equation}
$\Psi$ being the characteristic exponent of $\Gamma$, so that
$i\lambda c-\Psi(\lambda)=-\Psi_Y(\lambda)$. By \eqref{eq:distr2-gen},
\begin{equation}\label{eq:R-is-M}
R(\lambda)=\frac{\varphi_M(\lambda)}{\varphi_{-\overline H}(\lambda)} ,
\end{equation}
which is Lemma~\ref{lem:adjoint} read through transforms, and which reduces to
$R=\varphi_M$ exactly when $\overline H\equiv0$.
\end{remark}

\begin{remark}[Why an all-time supremum]\label{rem:PK}
Equation \eqref{eq:pow2-mid} is the balance equation of the tagged particle in
the wave frame, whose generator is \eqref{eq:tagged-generator}: the travelling
wave is the stationary law of a queue whose \emph{arrivals} are the
interactions, pushing the particle up by $X$, and whose \emph{server} is the
motion of the frame, pulling it down at speed $c$ against the drift
$\gamma_1$. The arrival rate $h(y)$ decreases with the content, so this is a
generalisation of M/G/1 rather than an instance of it.


The supremum is Reich's identity \cite{Reich58}. Solving the balance equation
means inverting the generator of the free motion, and for a driver with
negative drift that inverse is Lemma~\ref{lem:adjoint}, which for a spectrally
positive $Y$ is \eqref{eq:PK}. 

\end{remark}

With notation $\eta=-\chi$, \eqref{eq:distr2-gen} can be rewritten as
\begin{equation} \label{eq:distr3-gen}
\eta+\overline H \overset{d}{=} \max\{\eta_1,\eta_2\} - A,
\qquad A:=\overline{X}+M\ge0 ,
\end{equation}
which, when $\Gamma$ has no negative jumps and $\overline H\equiv0$, is
\begin{equation} \label{eq:distr3}
\eta \overset{d}{=} \max\{\eta_1,\eta_2\} - \overline{X}-M,
\end{equation}
a particular case of equation \eqref{eq:main_deq_general} with
$A=\overline{X}+M>0$; all the results of Section~\ref{sec:distr_equations} are
then applicable. For a two-sided driver the extra term $\overline H$ takes
\eqref{eq:distr3-gen} outside that framework, as we now explain.

\begin{remark}[Negative jumps]\label{rem:levy_case}
Proposition~\ref{prop:equiv_power2} holds for an arbitrary L\'evy driver, so
the decomposition is not what fails when $\Gamma$ has negative jumps; the
passage to Lemma~\ref{lemma:deq_main} is. Substituting the ansatz
$\Prob(\eta\le x)=\E[e^{-Ve^{-\gamma x}}]$ into \eqref{eq:distr3-gen} gives,
instead of \eqref{eq:deq_main_aux}, the modified fixed point
\begin{equation}\label{eq:modified-fp}
V\,e^{\gamma\overline H}\;\overset{d}{=}\;(V_1+V_2)\,e^{-\gamma A},
\qquad\text{with criticality}\qquad
\E\bigl[e^{-\gamma A}\bigr]=\tfrac12\,\E\bigl[e^{\gamma\overline H}\bigr],
\end{equation}
$\gamma$ being additionally required to lie below the abscissa of convergence
of $\overline H$. This is a smoothing transform carrying an independent random
factor on the \emph{left}, and the results of
\cite{liu1999asymptotic,Liu2001,Bigginsetal91} do not cover it. Studying
\eqref{eq:modified-fp} is the natural open problem left by this work.

\smallskip
\noindent\emph{Negative jumps flatten the wave.} The criticality relation in
\eqref{eq:modified-fp} is a moment identity, and is therefore available
without any existence theory. Since $\overline H\ge0$, its right-hand side is
at least $\tfrac12$ and non-decreasing in $\gamma$, whereas
$\gamma\mapsto\E[e^{-\gamma A}]$ decreases from $1$; the resulting $\gamma$ is
consequently strictly smaller than the root of $\E[e^{-\gamma A}]=\tfrac12$,
which is the value the same $A$ would produce for a spectrally positive
driver. Negative jumps in the driver thus lower the decay rate of the wave:
the travelling wave has a strictly heavier tail than the one-sided formula
predicts.

\smallskip
\noindent\emph{The transform reading.} The same phenomenon is visible in
\eqref{eq:R-is-M}: for a two-sided driver $R=\varphi_M/\varphi_{-\overline H}$
is a \emph{deconvolution}, hence the Fourier transform of a signed measure of
total mass one, and need not be a characteristic function.


\end{remark}

\subsection{Characterisation of the travelling wave}
Applying Lemma \ref{lemma:deq_main} to equation \eqref{eq:distr2}, we obtain the following characterisation:

\begin{theorem}[Travelling Wave Characterisation] \label{thm:TW_char}
For the power-of-2 model with a spectrally positive driver, so that
Corollary~\ref{cor:equiv_pow2_sp} applies, the travelling wave profile
$h(x) = \Prob(\chi > x)$ admits the representation
\begin{equation} \label{eq:TW_representation}
h(x) = \E\left[\exp\left(-V e^{\gamma x}\right)\right],
\end{equation}
where $\gamma > 0$ is a positive parameter which satisfies
\begin{equation} \label{eq:gamma_eq}
\frac{1}{2} = \E\left[\exp\left(-\gamma A\right)\right],
\end{equation}
and $V > 0$ is a random variable satisfying the linear equation
\begin{equation} \label{eq:V_equation}
V \overset{d}{=} (V_1 + V_2) \exp\left(-\gamma A\right),
\end{equation}
with $V_1, V_2$ independent copies of $V$, independent of $A = \overline{X} + M$,
where $M=\sup_{t\ge0}(\Gamma(t)-ct)$ is as in
Proposition~\ref{prop:equiv_power2}; recall that $A>0$ almost surely and has a
density, both because $\overline X$ does.
\end{theorem}

\begin{proof}
By Corollary~\ref{cor:equiv_pow2_sp} and \eqref{eq:distr3}, the variable
$\eta=-\chi$ solves \eqref{eq:main_deq_general} with $A=\overline X+M$, so
Lemma~\ref{lemma:deq_main} gives
$\Prob(\eta\le x)=\E\bigl[e^{-Ve^{-\gamma x}}\bigr]$; since
$h(x)=\Prob(\chi>x)=\Prob(\eta<-x)$, this is \eqref{eq:TW_representation}.

It remains to check the hypotheses of
Proposition~\ref{prop:linear-existence} for this $A$. The function
$\Phi(\gamma):=\E[e^{-\gamma A}]$ decreases strictly from $\Phi(0)=1$ to $0$,
because $A>0$ almost surely, so \eqref{eq:gamma_eq} has a unique solution
$\gamma>0$; and the supercritical condition $\E[Ae^{-\gamma A}]>0$ holds
automatically, again because $A>0$. Hence $V$ exists, is unique up to a
multiplicative constant and has finite mean; it is realised, as in
Remark~\ref{rem:BRW}, by the additive martingale of the branching random walk
with displacements distributed as $\gamma A$.
\end{proof}

\subsection{Properties of the Solution}

The characterisation in Theorem \ref{thm:TW_char} allows us to derive important properties of the travelling wave.

\begin{corollary}\label{cor:tail}
Under the hypotheses of Theorem~\ref{thm:TW_char}, write
$U:=\min\{\chi_1,\chi_2\}$, so that $\chi\overset{d}{=}U+A$ with $U$ and
$A=\overline X+M$ independent, and put
$\mu^*:=\sup\{\mu\ge0:\E[e^{\mu A}]<\infty\}$.
\begin{enumerate}
\item[(a)] $1-h(x)\sim\E[V]\,e^{\gamma x}$ as $x\to-\infty$.
\item[(b)] Assume $\mu^*>0$. Then $\E[e^{\beta\chi}]<\infty$ for $\beta<\mu^*$
and $\E[e^{\beta\chi}]=\infty$ for $\beta>\mu^*$, while $\E[e^{\beta U}]<\infty$
for every $\beta<2\mu^*$. If moreover
$\Prob(A>t)\sim L\,t^{k}e^{-\mu^*t}$ as $t\to\infty$, for some $L>0$ and
$k\ge0$, then
\[
h(x)\;\sim\;L\,\E\bigl[e^{\mu^*U}\bigr]\,x^{k}e^{-\mu^*x}
\qquad(x\to+\infty).
\]
\end{enumerate}
\end{corollary}

\begin{proof}
As for Lemma~\ref{lem:mu-star}, applied to $\chi$ in place of $\eta$: the two
statements are the same one read through the change of variable $\eta=-\chi$.
The moment hypothesis of Lemma~\ref{lem:mu-star}(b) is automatic here, since
$\E[V]<\infty$ by Theorem~\ref{thm:TW_char}.
\end{proof}

\begin{example}[Exponential jumps: the tail, explicitly]\label{ex:pow2-tail}
Take $X\sim\mathrm{Exp}(1)$ and $\sigma^2=1$. Then $\E[X]=1$, $c=1/2$,
$\overline X\sim\mathrm{Exp}(1)$ and, the driver being a driftless Brownian
motion of unit variance, $M\sim\mathrm{Exp}(1)$ by
Proposition~\ref{prop:equiv_power2}, so that
$A=\overline X+M\sim\Gamma(2,1)$ and \eqref{eq:gamma_eq} reads
\[
\E\bigl[e^{-\gamma A}\bigr]=\frac{1}{(1+\gamma)^2}=\frac12
\qquad\Longleftrightarrow\qquad
\gamma=\sqrt2-1,
\]
 Here
$\E[e^{\mu A}]=(1-\mu)^{-2}$, so $\mu^*=1>\gamma$, and Corollary~\ref{cor:tail} applies with $k=1$:
the decay rate is $\mu^*=1$, and since
$\Prob(A>t)=(1+t)e^{-t}$,
\[
h(x)=\E\bigl[(1+x-U)\,e^{-(x-U)}\bigr]\sim\E\bigl[e^{U}\bigr]\,x\,e^{-x}
\qquad(x\to+\infty),
\qquad U=\min\{\chi_1,\chi_2\},
\]
the expectation being finite because $U$ decays at rate $2$. The profile thus
decays like $xe^{-x}$. The factor $x$ comes from the $\Gamma(2,1)$ tail of
$A$, not from criticality: the power-of-2 wave is always supercritical
(Remark~\ref{rem:sign-of-A}).
\end{example}

\section{A model with mixed updates}\label{sec:mixed}

We close with a model in which the two mechanisms are combined inside a single
update, and which the method of Section~\ref{sec:power-of-2} solves in closed
form. Each particle carries two independent Poisson clocks, of rates
$\lambda_1$ and $\lambda_2$, and at every epoch draws an independent jump size
$J\ge0$ from a fixed law, normalised so that $\E[J]=1$.

\begin{itemize}
\item At a \emph{type-1} epoch the particle simply moves up by $J$.
\item At a \emph{type-2} epoch it contacts a particle chosen uniformly at
random and moves up by $J$, but not beyond the contacted particle: from state
$x$, with the contacted particle at $y$, the new state is
$\max\{x,\min\{x+J,y\}\}$, so that the increment is
$\min\{J,(y-x)^+\}$ and a particle contacting a lower one does not move.
\end{itemize}

Type-2 updates thus interpolate between the two mechanisms of this paper: the
jump is the power-of-2 push, the cap is the synchronisation rule. The model
describes, for instance, an industry in which type-1 is independent innovation
and type-2 is imitation of a competitor's technology, which cannot take a firm
beyond the firm it copies; or a peer-to-peer network in which type-1 is a
download from a complete central node and type-2 a download from a peer.

Writing $\lambda:=\lambda_1+\lambda_2$ for the total update rate, the
tagged particle initiates at rate $\lambda$ and
$\bar\theta(x,y)=I\,J+(1-I)\min\{J,(y-x)^+\}$ with $I$ Bernoulli of parameter
$\lambda_1/\lambda$, so Proposition~\ref{prop:TW-FP} with $\alpha=\lambda$ gives
\begin{equation}\label{eq:mixed-tagged}
\chi\overset{d}{=}\chi_1+I\,J+(1-I)\min\bigl\{J,(\chi_2-\chi_1)^+\bigr\}-cT,
\qquad T\sim\mathrm{Exp}(\lambda),
\end{equation}
all the variables on the right being independent and $\chi_1,\chi_2$ independent
copies of $\chi$. There is no diffusive component here, so the driver in the
wave frame is the deterministic drift $Y(t)=-ct$. The hypotheses of
Lemma~\ref{lem:bary-gen} hold as soon as $\E[J^2]<\infty$: the increment is
bounded by $J$, so that $\beta_2\le\E[J^2]$ and \eqref{eq:G} is satisfied,
while $\beta(x,y)=\lambda^{-1}\bigl(\lambda_1+\lambda_2\,\E\min\{J,(y-x)^+\}\bigr)$
is continuous.

\begin{proposition}[Speed and reduction to a max-type equation]
\label{prop:mixed}
Let $h(x)=\Prob(\chi>x)$. Any solution of \eqref{eq:mixed-tagged} has an
absolutely continuous law, and
\begin{equation}\label{eq:mixed-speed}
c=\frac{\lambda_2}{\nu},\qquad \nu:=\log\Bigl(1+\frac{\lambda_2}{\lambda_1}\Bigr),
\end{equation}
which does not depend on the law of $J$, and \eqref{eq:mixed-tagged} is then
equivalent to
\begin{equation}\label{eq:mixed-min}
\chi\overset{d}{=}\min_{1\le i\le K^+}\bigl(\chi_i+\overline J_i\bigr),
\end{equation}
where $\overline J_1,\overline J_2,\dots$ are independent copies of the
equilibrium jump $\overline J$, with density $\Prob(J>u)$ on $(0,\infty)$;
$\chi_1,\chi_2,\dots$ are independent copies of $\chi$; and $K^+$ is an
independent zero-truncated Poisson variable of parameter $\nu$,
\[
\Prob(K^+=k)=\frac{e^{-\nu}\nu^k}{k!\,(1-e^{-\nu})},\qquad k\ge1 .
\]
\end{proposition}

\begin{proof}
\emph{Step 1.} In the wave frame the tagged particle has generator
\[
\mathcal Lf(y)=-c\,f'(y)
+\lambda_1\,\E\bigl[f(y+J)-f(y)\bigr]
+\lambda_2\,\E\bigl[f\bigl(y+\min\{J,(\chi_2-y)^+\}\bigr)-f(y)\bigr],
\]
with $\chi_2\sim h$ independent, and, exactly as in Step~2 of the proof of
Proposition~\ref{prop:equiv_power2}, \eqref{eq:mixed-tagged} says that $h$ is
invariant for $\mathcal L$, i.e.\ $\E[\mathcal Lf(\chi)]=0$. Writing
$f(y+j)-f(y)=\int_0^{j}f'(y+u)\,du$ and using Fubini, the type-1 term is
$\int_0^\infty f'(y+u)\Prob(J>u)\,du=\E[f'(y+\overline J)]$, by
\eqref{eq:equilibrium-jump} and $\E[J]=1$. For the type-2 term the increment
exceeds $u$ exactly when $J>u$ and $\chi_2>y+u$, two independent events, so
that
\[
\E\bigl[f\bigl(y+\min\{J,(\chi_2-y)^+\}\bigr)\bigr]-f(y)
=\int_0^\infty f'(y+u)\,\Prob(J>u)\,h(y+u)\,du
=\E\bigl[f'(y+\overline J)\,h(y+\overline J)\bigr].
\]
Hence, for every $f\in\mathcal T$, with $g:=f'$,
\begin{equation}\label{eq:mixed-balance}
c\;\E\bigl[g(\chi)\bigr]
=\E\Bigl[g\bigl(\chi+\overline J\bigr)\,\bigl(\lambda_1+\lambda_2\,h(\chi+\overline J)\bigr)\Bigr],
\end{equation}
Both sides of \eqref{eq:mixed-balance} are integrals of $g$ against finite
measures. If $f\in\mathcal T$ has Fourier density $g_0$, then $g=f'$ has
Fourier density $i\lambda g_0(\lambda)$, so the Fourier transforms of the two
measures agree at every $\lambda\neq0$, hence everywhere by continuity. This
gives the identity between measures
\begin{equation}\label{eq:mixed-measure}
c\;\Prob(\chi\in dz)=\bigl(\lambda_1+\lambda_2h(z)\bigr)\,\Prob\bigl(\chi+\overline J\in dz\bigr).
\end{equation}

\emph{Step 2: absolute continuity, and the speed.} Integrating
\eqref{eq:mixed-measure} over $\mathbb R$ gives
$c=\lambda_1+\lambda_2\,\E[h(\chi+\overline J)]>0$. Since $\overline J$ has a
density, so does $\chi+\overline J$; and \eqref{eq:mixed-measure} exhibits the
law of $\chi$ as absolutely continuous with respect to that of
$\chi+\overline J$, hence as absolutely continuous. This proves the first
assertion, and in particular $h$ is continuous.

Read in the other direction, \eqref{eq:mixed-measure} says that the law of
$\chi+\overline J$ has density $c/(\lambda_1+\lambda_2h(z))$ with respect to
that of $\chi$; its total mass being $1$,
\[
c\;\E\Bigl[\frac{1}{\lambda_1+\lambda_2h(\chi)}\Bigr]=1 .
\]
As $h$ is a continuous survival function, $U:=h(\chi)$ is uniform on $(0,1)$,
so the expectation is a plain integral,
\[
c\int_0^1\frac{du}{\lambda_1+\lambda_2u}
=\frac{c}{\lambda_2}\log\frac{\lambda_1+\lambda_2}{\lambda_1}=1,
\]
which is \eqref{eq:mixed-speed}. The profile has disappeared from the
computation, which is why the speed does not depend on the law of $J$.

\emph{Step 3: the reduction.} Let $G(x):=\Prob(\chi+\overline J\le x)$.
Integrating \eqref{eq:mixed-measure} over $(-\infty,x]$ and substituting
$U=h(\chi)$ again, this time on the event $\{\chi\le x\}=\{U\ge h(x)\}$,
\[
G(x)=c\int_{h(x)}^1\frac{du}{\lambda_1+\lambda_2u}
=\frac1\nu\,\log\frac{\lambda_1+\lambda_2}{\lambda_1+\lambda_2h(x)} ,
\]
by \eqref{eq:mixed-speed}. Exponentiating, and using
$\lambda_1/\lambda=e^{-\nu}$ and $\lambda_2/\lambda=1-e^{-\nu}$,
\begin{equation}\label{eq:mixed-poisson}
e^{-\nu G(x)}=\frac{\lambda_1+\lambda_2h(x)}{\lambda}
=e^{-\nu}+\bigl(1-e^{-\nu}\bigr)h(x).
\end{equation}
Now let $K$ be Poisson of parameter $\nu$ and, conditionally on $K$, let
$Y_1,\dots,Y_K$ be independent with distribution function $G$, the minimum of
an empty family being $+\infty$. The left-hand side of
\eqref{eq:mixed-poisson} is the probability generating function of $K$
evaluated at $1-G(x)$, that is $\Prob(\min_{i\le K}Y_i>x)$, while the
right-hand side is $\Prob(K=0)+\Prob(K\ge1)h(x)$. Since the minimum exceeds
$x$ automatically on $\{K=0\}$, conditioning on $\{K\ge1\}$ leaves
$\Prob(\min_{i\le K^+}Y_i>x)=h(x)$, which is \eqref{eq:mixed-min} because
$Y_i\overset{d}{=}\chi_i+\overline J_i$. Every step is an equivalence.
\end{proof}

Equation \eqref{eq:mixed-min} is \eqref{eq:main_deq_general} with a random
number of branches: writing $\eta=-\chi$ it reads
\begin{equation}\label{eq:mixed-max}
\eta\overset{d}{=}\max_{1\le i\le K^+}\bigl(\eta_i-\overline J_i\bigr),
\end{equation}
and the proof of Lemma~\ref{lemma:deq_main} applies verbatim, the square of
the distribution function being replaced by its $K^+$-th moment: $\eta$ solves
\eqref{eq:mixed-max} with $\Prob(\eta\le x)=\E[e^{-Ve^{-\gamma x}}]$ if and
only if
\begin{equation}\label{eq:mixed-linear}
V\overset{d}{=}\sum_{i=1}^{K^+}e^{-\gamma\overline J_i}\,V_i .
\end{equation}
This is again a smoothing transform, now with $K^+$ summands, and the results
quoted in Proposition~\ref{prop:linear-existence} apply in that generality
\cite{Liu2001,liu1999asymptotic,Bigginsetal91}. Taking expectations in
\eqref{eq:mixed-linear}, and using
$\E[K^+]=\nu/(1-e^{-\nu})=\lambda/c$ together with
$\E[e^{-\gamma\overline J}]=\bigl(1-\E[e^{-\gamma J}]\bigr)/\gamma$, the
criticality condition $\E[K^+]\,\E[e^{-\gamma\overline J}]=1$ becomes
\begin{equation}\label{eq:mixed-gamma}
\gamma\,c=\lambda\bigl(1-\E[e^{-\gamma J}]\bigr),
\end{equation}
the analogue here of $D(\gamma)=2$. The right-hand side of
\eqref{eq:mixed-gamma} is concave and increasing in $\gamma$, with slope
$\lambda\E[J]=\lambda$ at the origin and limit $\lambda$ at infinity, while the
left-hand side is a line of slope $c<\lambda$; hence \eqref{eq:mixed-gamma}
has exactly one root $\gamma>0$. The profile is then
$h(x)=\E[e^{-Ve^{\gamma x}}]$, with $V$ the almost sure limit of the additive
martingale of the branching random walk whose offspring number is $K^+$ and
whose displacements are distributed as $\gamma\overline J$.

\begin{remark}[The two extremes]\label{rem:mixed-limits}
As $\lambda_2\to0$ one has $\nu\sim\lambda_2/\lambda_1$ and $c\to\lambda_1$,
the speed of a pure independent-jump system with $\E[J]=1$; correspondingly
$\E[K^+]\to1$ and $\gamma\to0$, the wave flattening out. As $\lambda_1\to0$,
$\nu\to\infty$ and $c\to0$: when every update is capped at another particle,
no particle can ever pass the current maximum, the front stalls, and the
speed \eqref{eq:mixed-speed} records this. The interesting regime is
therefore the interior one, in which an arbitrarily small rate of independent
innovation produces a strictly positive speed
$\lambda_2/\log(1+\lambda_2/\lambda_1)$, vanishing only logarithmically as
$\lambda_1\downarrow0$.
\end{remark}

\section{Concluding remarks and open problems}

We have characterised travelling waves of mean-field equations through
a distributional equation for a tagged particle seen in the frame of the wave,
solved by the correspondence of Lemma~\ref{lemma:deq_main} with a smoothing
transform. For the synchronisation mechanism this recovers the F-KPP waves
and, for compound Poisson motion, gives the profile for every jump law with an
exponential moment. For a model in which each update mixes the two mechanisms,
it produces the speed and the profile in closed form, through a max-type
equation with a random number of branches (Section~\ref{sec:mixed}). For the power-of-2
mechanism it identifies the wave as the stationary law of a queue: the balance
equation of the tagged particle yields the decomposition
\eqref{eq:distr2-gen} for an arbitrary L\'evy driver, exhibits the speed as
the stability condition of that queue, and gives the exact decay rate of each
tail.

Three questions are left open. When the driver has negative jumps, the
decomposition still holds but the fixed point acquires an extra random factor,
\eqref{eq:modified-fp}, which falls outside the class of smoothing transforms;
all we know is that it forces a strictly heavier tail than the one-sided
formula predicts. The sharp constant in the left tail of
Section~\ref{sec:prob-FKPP} is established only for Brownian synchronisation,
the identification of the exponent being general
(Remark~\ref{rem:left-tail}). Finally, Lemma~\ref{lemma:deq_main} characterises
the solutions of \eqref{eq:main_deq_general} of Gumbel-mixture form; that the
non-linear equation admits no others should follow from \cite{ABM12}, as
explained in Remark~\ref{rem:deq_main}, but a careful verification of their
hypotheses in our setting remains to be done.

\appendix

\section{The equilibrium ladder height and the adjoint relation}\label{app:ladder}

Let $Y$ be a L\'evy process with $\E|Y(1)|<\infty$ and
$\E\,Y(1)=-\mu<0$, so that $M:=\sup_{t\ge0}Y(t)$ is almost surely finite, and
let $\hat H$ be the descending ladder height process of $Y$, a subordinator
with Laplace exponent $\hat\varkappa$, drift $d$, L\'evy measure $\hat\Pi$ and
mean $m:=\E\,\hat H_1$. We write $\overline H\ge0$ for the associated
\emph{equilibrium}, or stationary overshoot, law,
\begin{equation}\label{eq:Hbar}
\Prob\bigl(\overline H\in du\bigr)
=\frac{1}{m}\Bigl(d\,\delta_0(du)+\hat\Pi\bigl((u,\infty)\bigr)\,du\Bigr),
\qquad\text{equivalently}\qquad
\E\bigl[e^{-s\overline H}\bigr]=\frac{\hat\varkappa(s)}{m\,s},
\end{equation}
a probability law which does not depend on the normalisation of the
Wiener--Hopf factorisation. In particular $\overline H\equiv0$ if and only if $\hat\Pi=0$, that is, if and
only if $Y$ has no negative jumps.

It is convenient to fix once and for all the class of test functions
\begin{equation}\label{eq:testclass}
\mathcal T:=\Bigl\{f(x)=\tfrac1{2\pi}\int_{\mathbb R}g(\lambda)e^{i\lambda x}\,d\lambda
\;:\;\int_{\mathbb R}(1+\lambda^2)\,|g(\lambda)|\,d\lambda<\infty\Bigr\}.
\end{equation}
Every $f\in\mathcal T$ is bounded and twice continuously differentiable with
bounded derivatives, so that $\mathcal A_Yf$ is defined and bounded;
$\mathcal T$ contains the Schwartz class and therefore separates finite signed
measures; and it is stable under the two operations used in
Section~\ref{sec:power-of-2}, namely translation and averaging against an
independent random variable, which act on $g$ by multiplication by
$e^{i\lambda a}$ and by a characteristic function respectively --- note that
the Schwartz class itself is \emph{not} stable under the second.

\begin{lemma}[Adjoint relation for the all-time supremum]\label{lem:adjoint}
In the setting above, for every $f\in\mathcal T$,
\begin{equation}\label{eq:adjoint}
\E\bigl[\mathcal A_Yf(M)\bigr]=-\mu\,\E\bigl[f'(-\overline H)\bigr],
\end{equation}
where $\mathcal A_Y$ is the generator of $Y$ and $\overline H$ is independent
of everything else. When $Y$ has no negative jumps this reads
$\E[\mathcal A_Yf(M)]=-\mu\,f'(0)$.
\end{lemma}

\begin{proof}
Since $|\Psi_Y(\lambda)|\le C(1+\lambda^2)$ for some $C<\infty$, both
$g\Psi_Y$ and $\lambda g$ are integrable, so Fubini gives
\[
\E\bigl[\mathcal A_Yf(M)\bigr]
=\frac1{2\pi}\int g(\lambda)\,\Psi_Y(\lambda)\,\E\bigl[e^{i\lambda M}\bigr]d\lambda,
\qquad
-\mu\,\E\bigl[f'(-\overline H)\bigr]
=\frac1{2\pi}\int g(\lambda)\,(-i\lambda\mu)\,\E\bigl[e^{-i\lambda\overline H}\bigr]d\lambda,
\]
so that it is enough to check that the integrands agree for every $\lambda$,
i.e.\ that
$\Psi_Y(\lambda)\E[e^{i\lambda M}]=-i\lambda\mu\,\E[e^{-i\lambda\overline H}]$.
By the Wiener--Hopf factorisation \cite{bertoin1996},
$-\Psi_Y(\lambda)=\varkappa(-i\lambda)\hat\varkappa(i\lambda)$ and
$\E[e^{i\lambda M}]=\varkappa(0)/\varkappa(-i\lambda)$, $\varkappa$ denoting the
ascending ladder exponent, while \eqref{eq:Hbar} gives
$\E[e^{-i\lambda\overline H}]=\hat\varkappa(i\lambda)/(m\,i\lambda)$. Hence
\begin{equation*}
\frac{-i\lambda\mu\,\E[e^{-i\lambda\overline H}]}{\Psi_Y(\lambda)}
=\frac{i\lambda\mu}{\varkappa(-i\lambda)\,\hat\varkappa(i\lambda)}
\cdot\frac{\hat\varkappa(i\lambda)}{m\,i\lambda}
=\frac{\mu}{m\,\varkappa(-i\lambda)},
\end{equation*}
which is $\E[e^{i\lambda M}]$ provided $\varkappa(0)=\mu/m$. That identity
follows by letting $\lambda\to0$ in
$-\Psi_Y(\lambda)=\varkappa(-i\lambda)\hat\varkappa(i\lambda)$, since
$-\Psi_Y(\lambda)\sim i\lambda\mu$ and $\hat\varkappa(i\lambda)\sim m\,i\lambda$.
Finally $m<\infty$, because $m=\mu/\varkappa(0)$ and $\varkappa(0)>0$, the ascending
ladder height being killed as $Y$ drifts to $-\infty$.
\end{proof}

\bibliographystyle{abbrv}
\bibliography{biblio}

@article{ABM12,
  author  = {G. Alsmeyer and J. D. Biggins and M. Meiners},
  title   = {The functional equation of the smoothing transform},
  journal = {Ann. Probab.},
  volume  = {40},
  number  = {5},
  pages   = {2069--2105},
  year    = {2012}
}

@book{Sato99,
  author    = {K. Sato},
  title     = {L\'evy Processes and Infinitely Divisible Distributions},
  series    = {Cambridge Studies in Advanced Mathematics},
  volume    = {68},
  publisher = {Cambridge University Press},
  year      = {1999}
}

@book{Asmussen03,
  author    = {Asmussen, S{\o}ren},
  title     = {Applied Probability and Queues},
  edition   = {2nd},
  series    = {Applications of Mathematics (New York)},
  volume    = {51},
  publisher = {Springer},
  address   = {New York},
  year      = {2003},
  note      = {Stochastic Modelling and Applied Probability},
  isbn      = {0-387-00211-1}
}

@article{Reich58,
  author  = {Reich, Edgar},
  title   = {On the integrodifferential equation of {T}ak\'acs {I}},
  journal = {Ann. Math. Statist.},
  volume  = {29},
  number  = {2},
  pages   = {563--570},
  year    = {1958},
  doi     = {10.1214/aoms/1177706632}
}

@incollection{Sznitman91,
  author    = {Sznitman, Alain-Sol},
  title     = {Topics in propagation of chaos},
  booktitle = {\'Ecole d'\'Et\'e de Probabilit\'es de Saint-Flour XIX --- 1989},
  editor    = {Hennequin, Paul-Louis},
  series    = {Lecture Notes in Mathematics},
  volume    = {1464},
  pages     = {165--251},
  publisher = {Springer},
  address   = {Berlin},
  year      = {1991},
  doi       = {10.1007/BFb0085169}
}

@book{bertoin1996,
  author = {Bertoin, Jean},
  title = {L\'{e}vy Processes},
  series = {Cambridge Tracts in Mathematics},
  volume = {121},
  publisher = {Cambridge University Press},
  address = {Cambridge},
  year = {1996}
}

@article{aldous2005,
  author = {Aldous, D. and Bandyopadhyay, A.},
  title = {A survey of max-type recursive distributional equations},
  journal = {Ann. Appl. Probab.},
  volume = {15},
  number = {2},
  pages = {1047--1110},
  year = {2005},
  doi = {10.1214/105051605000000142}
}

@article{Biggins1,
  author = {J. D. Biggins},
  title = {The growth and spread of the general branching random walk},
  journal = {Ann. Appl. Probab.},
  volume = {5},
  number = {4},
  pages = {1008--1024},
  year = {1995}
}

@incollection{Biggins2,
  author = {J. D. Biggins},
  title = {How fast does a general branching random walk spread?},
  booktitle = {Classical and modern branching processes (Minneapolis, MN, 1994)},
  volume = {84},
  series = {IMA Vol. Math. Appl.},
  pages = {19--39},
  publisher = {Springer, New York},
  year = {1997}
}

@article{Bigginsetal91,
  author = {J. D. Biggins and Boris D. Lubachevsky and Adam Shwartz and Alan Weiss},
  title = {A branching random walk with a barrier},
  journal = {Ann. Appl. Probab.},
  volume = {1},
  number = {4},
  pages = {573--581},
  year = {1991}
}

@article{B2,
  author = {M. Bramson},
  title = {Convergence of solutions of the Kolmogorov equation to travelling waves},
  journal = {Mem. Amer. Math. Soc.},
  volume = {44},
  number = {285},
  pages = {iv+190},
  year = {1983}
}

@article{BD,
  author = {{\'E}. Brunet and B. Derrida},
  title = {Effect of microscopic noise on front propagation},
  journal = {J. Statist. Phys.},
  volume = {103},
  number = {1-2},
  pages = {269--282},
  year = {2001}
}

@article{GJ2018,
  title={Front propagation and quasi-stationary distributions for one-dimensional L{\'e}vy processes},
  author={Groisman, Pablo and Jonckheere, Matthieu},
  journal={Electronic Communications in Probability},
  volume={23},
  pages={1--11},
  year={2018},
  publisher={Institute of Mathematical Statistics and Bernoulli Society}
}

@article{GJM,
  title={F-KPP Scaling limit and selection principle for a Brunet-Derrida type particle system},
  author={Groisman, Pablo and Jonckheere, Matthieu and Mart{\'i}nez, Juli{\'a}n},
  journal={ALEA: Latin American Journal of Probability and Mathematical Statistics},
  volume={17},
  number={1},
  pages={589},
  year={2020}
}

@article{Liu2001,
  author  = {Q. Liu},
  title   = {Asymptotic properties and absolute continuity of laws stable by random weighted mean},
  journal = {Stochastic Process. Appl.},
  volume  = {95},
  number  = {1},
  pages   = {83--107},
  year    = {2001}
}

@article{Fisher,
  author = {R. A. Fisher},
  title = {The wave of advance of advantageous genes},
  journal = {Ann. Eugenics},
  volume = {7},
  pages = {353--369},
  year = {1937}
}

@article{KPP,
  author = {A. Kolmogorov and I. Petrovsky and N. Piscounov},
  title = {Etude de l'equation de la diffusion avec croissance de la quantite de matiere et son application a un probleme biologique},
  journal = {Bull. Univ. Etat Moscou},
  volume = {A 1},
  number = {1},
  pages = {25},
  year = {1937}
}

@article{Ky,
  author = {A. E. Kyprianou},
  title = {A note on branching Levy processes},
  journal = {Stochastic Process. Appl.},
  volume = {82},
  number = {1},
  pages = {1--14},
  year = {1999}
}

@article{liu1999asymptotic,
  title={Asymptotic properties of supercritical age-dependent branching processes and homogeneous branching random walks},
  author={Liu, Quansheng},
  journal={Stochastic Processes and their Applications},
  volume={82},
  number={1},
  pages={61--87},
  year={1999},
  publisher={Elsevier}
}

@article{MK,
  author = {H. P. McKean},
  title = {Application of Brownian motion to the equation of Kolmogorov-Petrovskii-Piskunov},
  journal = {Comm. Pure Appl. Math.},
  volume = {28},
  number = {3},
  pages = {323--331},
  year = {1975}
}

@article{baryshnikov2025large,
  title={A large-scale particle system with independent jumps and distributed synchronization},
  author={Baryshnikov, Yuliy and Stolyar, Alexander},
  journal={Advances in Applied Probability},
  volume={57},
  number={2},
  pages={677--707},
  year={2025},
  publisher={Cambridge University Press}
}

@article{deMasi,
  title = {Hydrodynamics of the {$N$}-{BBM} process},
  author = {De Masi, A. and Ferrari, P.A. and Presutti, E. and Soprano-Loto, N.},
  year = {2017},
  journal = {arXiv:1707.00799}
}

@book{Liggett,
  author = {Liggett, T.M.},
  title = {Stochastic interacting systems: contact, voter and exclusion processes},
  publisher = {Springer-Verlag},
  series = {Grundlehren der Mathematischen Wissenschaften [Fundamental Principles of Mathematical Sciences]},
  volume = {324},
  year = {1999}
}

@book{E,
  author = {Etheridge, A.M.},
  title = {An introduction to superprocesses},
  publisher = {American Mathematical Society},
  year = {2000}
}

@article{DR,
  author = {Durrett, R. and Remenik, D.},
  title = {Brunet-derrida particle systems, free boundary problems and Wiener-Hopf equations},
  journal = {Ann. Probab.},
  number = {6},
  pages = {2043--2078},
  volume = {39},
  year = {2011}
}

@article{BD2,
  author = {Brunet, E. and Derrida, B.},
  title = {Shift in the velocity of a front due to a cutoff},
  journal = {Phys. Rev. E (3)},
  volume = {56},
  number = {3, part A},
  pages = {2597--2604},
  year = {1997}
}

@article{BDMM,
  author = {Brunet, E. and Derrida, B. and Mueller, A. H. and Munier, S.},
  title = {Noisy traveling waves: effect of selection on genealogies},
  journal = {Europhys. Lett.},
  volume = {76},
  number = {1},
  pages = {1--7},
  year = {2006}
}

@article{BDMM2,
  author = {Brunet, {\'E}. and Derrida, B. and Mueller, A. H. and Munier, S.},
  title = {Effect of selection on ancestry: an exactly soluble case and its phenomenological generalization},
  journal = {Phys. Rev. E (3)},
  volume = {76},
  number = {4},
  pages = {041104, 20},
  year = {2007}
}

@article{K,
  author = {Kolmogorov, AN},
  title = {On the solution of a biological problem},
  journal = {Proceedings of Tomsk University},
  pages = {7--12},
  volume = {2},
  year = {1938}
}

@article{S,
  author = {Sevastyanov, A.B.},
  title = {The theory of branching processes},
  journal = {Uspekhi Mat.},
  number = {6},
  pages = {47--99},
  volume = {6},
  year = {1951}
}

@book{B,
  author = {Billingsley, P.},
  title = {Convergence of probability measures},
  edition = {Second},
  publisher = {John Wiley \& Sons Inc.},
  series = {Wiley Series in Probability and Statistics: Probability and Statistics},
  year = {1999}
}

@article{M,
  author = {Maillard, Pascal},
  title = {Speed and fluctuations of {$N$}-particle branching {B}rownian motion with spatial selection},
  journal = {Probab. Theory Related Fields},
  volume = {166},
  number = {3-4},
  pages = {1061--1173},
  year = {2016}
}

@article{MMQ2,
  author = {Mueller, Carl and Mytnik, Leonid and Quastel, Jeremy},
  title = {Effect of noise on front propagation in reaction-diffusion equations of {KPP} type},
  journal = {Invent. Math.},
  volume = {184},
  number = {2},
  pages = {405--453},
  year = {2011}
}

@article{MMQ,
  author = {Mueller, C. and Mytnik, L. and Quastel, J.},
  title = {Small noise asymptotics of traveling waves},
  journal = {Markov Process. Related Fields},
  volume = {14},
  number = {3},
  pages = {333--342},
  year = {2008}
}

@article{DP,
  author = {Dumortier, F. and Popovi\'{c}, N. and Kaper, T.J.},
  title = {The critical wave speed for the {F}isher-{K}olmogorov-{P}etrowskii-{P}iscounov equation with cut-off},
  journal = {Nonlinearity},
  volume = {20},
  number = {4},
  pages = {855--877},
  year = {2007}
}

@article{BenDep,
  author = {Benguria, RD and Depassier, MC},
  title = {Speed of pulled fronts with a cutoff},
  journal = {Physical Review E},
  volume = {75},
  number = {5},
  pages = {051106},
  year = {2007},
  publisher = {APS}
}

@article{AW,
  author = {Aronson, D. G. and Weinberger, H. F.},
  title = {Multidimensional nonlinear diffusion arising in population genetics},
  journal = {Adv. in Math.},
  volume = {30},
  number = {1},
  pages = {33--76},
  year = {1978}
}

@article{U,
  author = {Uchiyama, K.},
  title = {The behavior of solutions of some nonlinear diffusion equations for large time},
  journal = {J. Math. Kyoto Univ.},
  volume = {18},
  number = {3},
  pages = {453--508},
  year = {1978}
}

@article{Bramson78,
  author = {Bramson, M.D.},
  title = {Maximal displacement of branching {B}rownian motion},
  journal = {Comm. Pure Appl. Math.},
  volume = {31},
  number = {5},
  pages = {531--581},
  year = {1978}
}

@article{Stolyar2023,
    author = {Stolyar, Alexander L.},
    title = {A Particle System with Mean-Field Interaction: Large-Scale Limit of Stationary Distributions},
    journal = {Stochastic Systems},
    year = {2023},
    volume = {13},
    number = {3},
    pages = {343--359},
    doi = {10.1287/stsy.2023.0108}
}

@article{Stolyar2023TW,
    author = {Stolyar, Alexander L.},
    title = {Large-scale behavior of a particle system with mean-field interaction: Traveling wave solutions},
    journal = {Advances in Applied Probability},
    year = {2023},
    volume = {55},
    number = {1},
    pages = {245--274},
    doi = {10.1017/apr.2022.35}
}

@article{Greenberg1995,
    author = {Greenberg, A. and Malyshev, V. and Popov, S.},
    title = {Stochastic model of massively parallel computation},
    journal = {Markov Processes and Related Fields},
    year = {1995},
    volume = {1},
    number = {4},
    pages = {473--490}
}

@inproceedings{Greenberg1996,
    author = {Greenberg, A. and Shenker, S. and Stolyar, A.},
    title = {Asynchronous updates in large parallel systems},
    booktitle = {Proceedings of the ACM SIGMETRICS International Conference on Measurement and Modeling of Computer Systems},
    year = {1996},
    pages = {91--103},
    publisher = {ACM}
}

@article{Balazs2014,
    author = {Balazs, M. and Racz, M. Z. and Toth, B.},
    title = {Modeling flocks and prices: Jumping particles with an attractive interaction},
    journal = {Annales de l'Institut Henri Poincaré, Probabilités et Statistiques},
    year = {2014},
    volume = {50},
    number = {2},
    pages = {425--454},
    doi = {10.1214/12-AIHP529}
}

@article{Hongler2015,
    author = {Hongler, M. O.},
    title = {Exact soliton-like probability measures for interacting jump processes},
    journal = {Mathematical Sciences},
    year = {2015},
    volume = {40},
    number = {1},
    pages = {62--66},
    doi = {10.1007/s40096-015-0156-5}
}

@article{Hongler2019,
    author = {Hongler, M. O. and Filliger, R.},
    title = {When do redundant requests reduce latency?},
    journal = {Methodology and Computing in Applied Probability},
    year = {2019},
    volume = {21},
    pages = {753--764},
    doi = {10.1007/s11009-018-9633-4}
}

\end{document}